\documentclass[12pt]{article}
\usepackage{theorem}
\usepackage{amssymb}
\usepackage{graphicx}
\usepackage[mathscr]{eucal}
\usepackage{amsbsy}
\usepackage{amsmath}
\usepackage{stmaryrd}
\usepackage{subfigure}
\newtheorem{prop}{}[section]

{\theorembodyfont{\upshape} \newtheorem{rema}[prop]{}}
\newcommand{\boma}[1]{{\mbox{\boldmath $#1$} }}

\begin{document}
\newcommand{\uper}[1]{\stackrel{\barray{c} {~} \\ \mbox{\footnotesize{#1}}\farray}{\longrightarrow} }
\newcommand{\nop}[1]{ \|#1\|_{\piu} }
\newcommand{\no}[1]{ \|#1\| }
\newcommand{\nom}[1]{ \|#1\|_{\meno} }
\newcommand{\UU}[1]{e^{#1 \AA}}
\newcommand{\UD}[1]{e^{#1 \Delta}}
\newcommand{\bb}[1]{\mathbb{{#1}}}
\newcommand{\HO}[1]{\bb{H}^{{#1}}}
\newcommand{\Hz}[1]{\bb{H}^{{#1}}_{\zz}}
\newcommand{\Hs}[1]{\bb{H}^{{#1}}_{\ss}}
\newcommand{\Hg}[1]{\bb{H}^{{#1}}_{\gg}}
\newcommand{\HM}[1]{\bb{H}^{{#1}}_{\so}}
\newcommand{\Hm}[1]{H^{{#1}}}
\def\ohz{O_h \ltimes \Tt}
\def\ohperz{O_h \times \Tt}
\def\t{{\tt{T}}}
\def\Isot{{\mathcal H}}
\def\PIsot{{\mathcal H}^{-}}
\def\Isott{{\mathcal H}_R}
\def\PIsott{{\mathcal H}^{-}_R}
\def\uno{{\bf 1}}
\def\to{t_{\circ}}
\def\to{T_{\circ}}
\def\ts{T_{*}}
\def\es{S_{*}}
\def\Ti{\mathscr{T}}
\def\Pade{Pad\'e~}
\def\bbb{\interleave}
\def\Aa{\mathscr{A}}
\def\ef{\psi}
\def\fun{\mathcal{F}}
\def\fun{{\tt f}}
\def\tvainf{\vspace{-0.4cm} \barray{ccc} \vspace{-0,1cm}{~}
\\ \vspace{-0.2cm} \longrightarrow \\ \vspace{-0.2cm} \scriptstyle{T \vain + \infty} \farray}
\def\De{F}
\def\er{\epsilon}
\def\erd{\er_0}
\def\Tn{T_{\star}}
\def\Tc{T_{\tt{c}}}
\def\Tb{T_{\tt{b}}}
\def\Tl{\mathscr{T}}
\def\Tm{T}
\def\Ta{T_{\tt{a}}}
\def\ua{u_{\tt{a}}}
\def\Tg{T_{G}}
\def\Tgg{T_{I}}
\def\Tw{T_{w}}
\def\Ts{T_{\Ss}}
\def\Ss{T_{\Ss}}
\def\Tr{\Tl}
\def\Sp{\Ss'}
\def\Tsp{T_{\Sp}}
\def\vsm{\vspace{-0.1cm}\noindent}
\def\comple{\scriptscriptstyle{\complessi}}
\def\ug{u_G}
\def\nume{0.407}
\def\numerob{0.00724}
\def\deln{7/10}
\def\delnn{\dd{7 \over 10}}
\def\e{c}
\def\p{p}
\def\z{z}
\def\symd{{\mathfrak S}_d}
\def\del{\omega}
\def\Del{\delta}
\def\Di{\Delta}
\def\Ss{{\mathscr{S}}}
\def\Ww{{\mathscr{W}}}
\def\mmu{\hat{\mu}}
\def\rot{\mbox{rot}\,}
\def\curl{\mbox{curl}\,}
\def\Mm{\mathscr M}
\def\XS{\boma{x}}
\def\TS{\boma{t}}
\def\Lam{\boma{\eta}}
\def\DS{\boma{\rho}}
\def\KS{\boma{k}}
\def\LS{\boma{\lambda}}
\def\PR{\boma{p}}
\def\VS{\boma{v}}
\def\ski{\! \! \! \! \! \! \! \! \! \! \! \! \! \!}
\def\h{L}
\def\EM{M}
\def\EMP{M'}
\def\R{R}
\def\Rr{{\mathscr{R}}}
\def\Zz{{\mathscr{Z}}}
\def\FFf{\mathscr{F}}
\def\A{F}
\def\Xim{\Xi_{\meno}}
\def\Ximn{\Xi_{n-1}}
\def\lan{\lambda}
\def\om{\omega}
\def\Om{\Omega}
\def\Sim{\Sigm}
\def\Sip{\Delta \Sigm}
\def\Sigm{{\mathscr{S}}}
\def\Ki{{\mathscr{K}}}
\def\Hi{{\mathscr{H}}}
\def\zz{{\scriptscriptstyle{0}}}
\def\ss{{\scriptscriptstyle{\Sigma}}}
\def\gg{{\scriptscriptstyle{\Gamma}}}
\def\so{\ss \zz}
\def\Dv{\bb{\DD}'}
\def\Dz{\bb{\DD}'_{\zz}}
\def\Ds{\bb{\DD}'_{\ss}}
\def\Dsz{\bb{\DD}'_{\so}}
\def\Dg{\bb{\DD}'_{\gg}}
\def\Ls{\bb{L}^2_{\ss}}
\def\Lg{\bb{L}^2_{\gg}}
\def\bF{{\bb{V}}}
\def\Fz{\bF_{\zz}}
\def\Fs{\bF_\ss}
\def\Fg{\bF_\gg}
\def\Pre{P}
\def\UUU{{\mathcal U}}
\def\fiapp{\phi}
\def\PU{P1}
\def\PD{P2}
\def\PT{P3}
\def\PQ{P4}
\def\PC{P5}
\def\PS{P6}
\def\Q{P6}
\def\X{Q2}
\def\Xp{Q3}
\def\Vi{V}
\def\bVi{\bb{V}}
\def\K{V}
\def\Ks{\bb{\K}_\ss}
\def\Kz{\bb{\K}_0}
\def\KM{\bb{\K}_{\, \so}}
\def\HGG{\bb{H}^\G}
\def\HG{\bb{H}^\G_{\so}}
\def\EG{{\mathfrak{P}}^{\G}}
\def\G{G}
\def\de{\delta}
\def\esp{\sigma}
\def\dd{\displaystyle}
\def\LP{\mathfrak{L}}
\def\dive{\mbox{div}}
\def\la{\langle}
\def\ra{\rangle}
\def\um{u_{\meno}}
\def\uv{\mu_{\meno}}
\def\Fp{ {\textbf F_{\piu}} }
\def\Ff{ {\textbf F} }
\def\Fm{ {\textbf F_{\meno}} }
\def\Eb{ {\textbf E} }
\def\piu{\scriptscriptstyle{+}}
\def\meno{\scriptscriptstyle{-}}
\def\omeno{\scriptscriptstyle{\ominus}}
\def\Tt{ {\mathscr T} }
\def\Xx{ {\textbf X} }
\def\Yy{ {\textbf Y} }
\def\Ee{ {\textbf E} }
\def\VP{{\mbox{\tt VP}}}
\def\CP{{\mbox{\tt CP}}}
\def\cp{$\CP(f_0, t_0)\,$}
\def\cop{$\CP(f_0)\,$}
\def\copn{$\CP_n(f_0)\,$}
\def\vp{$\VP(f_0, t_0)\,$}
\def\vop{$\VP(f_0)\,$}
\def\vopn{$\VP_n(f_0)\,$}
\def\vopdue{$\VP_2(f_0)\,$}
\def\leqs{\leqslant}
\def\geqs{\geqslant}
\def\mat{{\frak g}}
\def\tG{t_{\scriptscriptstyle{G}}}
\def\tN{t_{\scriptscriptstyle{N}}}
\def\TK{t_{\scriptscriptstyle{K}}}
\def\CK{C_{\scriptscriptstyle{K}}}
\def\CN{C_{\scriptscriptstyle{N}}}
\def\CG{C_{\scriptscriptstyle{G}}}
\def\CCG{{\mathscr{C}}_{\scriptscriptstyle{G}}}
\def\tf{{\tt f}}
\def\ti{{\tt t}}
\def\ta{{\tt a}}
\def\tc{{\tt c}}
\def\tF{{\tt R}}
\def\C{{\mathscr C}}
\def\P{{\mathscr P}}
\def\V{{\mathscr V}}
\def\TI{\tilde{I}}
\def\TJ{\tilde{J}}
\def\Lin{\mbox{Lin}}
\def\Hinfc{ H^{\infty}(\reali^d, \complessi) }
\def\Hnc{ H^{n}(\reali^d, \complessi) }
\def\Hmc{ H^{m}(\reali^d, \complessi) }
\def\Hac{ H^{a}(\reali^d, \complessi) }
\def\Dc{\DD(\reali^d, \complessi)}
\def\Dpc{\DD'(\reali^d, \complessi)}
\def\Sc{\SS(\reali^d, \complessi)}
\def\Spc{\SS'(\reali^d, \complessi)}
\def\Ldc{L^{2}(\reali^d, \complessi)}
\def\Lpc{L^{p}(\reali^d, \complessi)}
\def\Lqc{L^{q}(\reali^d, \complessi)}
\def\Lrc{L^{r}(\reali^d, \complessi)}
\def\Hinfr{ H^{\infty}(\reali^d, \reali) }
\def\Hnr{ H^{n}(\reali^d, \reali) }
\def\Hmr{ H^{m}(\reali^d, \reali) }
\def\Har{ H^{a}(\reali^d, \reali) }
\def\Dr{\DD(\reali^d, \reali)}
\def\Dpr{\DD'(\reali^d, \reali)}
\def\Sr{\SS(\reali^d, \reali)}
\def\Spr{\SS'(\reali^d, \reali)}
\def\Ldr{L^{2}(\reali^d, \reali)}
\def\Hinfk{ H^{\infty}(\reali^d, \KKK) }
\def\Hnk{ H^{n}(\reali^d, \KKK) }
\def\Hmk{ H^{m}(\reali^d, \KKK) }
\def\Hak{ H^{a}(\reali^d, \KKK) }
\def\Dk{\DD(\reali^d, \KKK)}
\def\Dpk{\DD'(\reali^d, \KKK)}
\def\Sk{\SS(\reali^d, \KKK)}
\def\Spk{\SS'(\reali^d, \KKK)}
\def\Ldk{L^{2}(\reali^d, \KKK)}
\def\Knb{K^{best}_n}
\def\sc{\cdot}
\def\k{\mbox{{\tt k}}}
\def\x{\mbox{{\tt x}}}
\def\g{ {\textbf g} }
\def\QQQ{ {\textbf Q} }
\def\AAA{ {\textbf A} }
\def\gr{\mbox{gr}}
\def\sgr{\mbox{sgr}}
\def\loc{\mbox{loc}}
\def\PZ{{\Lambda}}
\def\PZAL{\mbox{P}^{0}_\alpha}
\def\epsilona{\epsilon^{\scriptscriptstyle{<}}}
\def\epsilonb{\epsilon^{\scriptscriptstyle{>}}}
\def\lgraffa{ \mbox{\Large $\{$ } \hskip -0.2cm}
\def\rgraffa{ \mbox{\Large $\}$ } }
\def\restriction{\upharpoonright}
\def\M{{\scriptscriptstyle{M}}}
\def\m{m}
\def\Fre{Fr\'echet~}
\def\I{{\mathcal N}}
\def\ap{{\scriptscriptstyle{ap}}}
\def\fiap{\varphi_{\ap}}
\def\dfiap{{\dot \varphi}_{\ap}}
\def\DDD{ {\mathfrak D} }
\def\BBB{ {\textbf B} }
\def\EEE{ {\textbf E} }
\def\GGG{ {\textbf G} }
\def\TTT{ {\textbf T} }
\def\KKK{ {\textbf K} }
\def\HHH{ {\textbf K} }
\def\FFi{ {\bf \Phi} }
\def\GGam{ {\bf \Gamma} }
\def\sc{ {\scriptstyle{\bullet} }}
\def\a{a}
\def\ep{\epsilon}
\def\c{\kappa}
\def\parn{\par \noindent}
\def\teta{M}
\def\elle{L}
\def\ro{\rho}
\def\al{\alpha}
\def\si{\sigma}
\def\be{\beta}
\def\ga{\gamma}
\def\te{\vartheta}
\def\ch{\chi}
\def\et{\eta}
\def\complessi{{\bf C}}
\def\len{{\bf L}}
\def\reali{{\bf R}}
\def\interi{{\bf Z}}
\def\razionali{{\bf Q}}
\def\Z{{\bf Z}}
\def\naturali{{\bf N}}
\def\Sfe{ {\bf S} }
\def\To{ {\bf T} }
\def\Td{ {\To}^d }
\def\Tt{ {\To}^3 }
\def\Zd{ \interi^d }
\def\Zt{ \interi^3 }
\def\Zet{{\mathscr{Z}}}
\def\Ze{\Zet^d}
\def\T1{{\textbf To}^{1}}
\def\ee{{E}}
\def\FF{\mathcal F}
\def\FFu{ {\textbf F_{1}} }
\def\FFd{ {\textbf F_{2}} }
\def\GG{{\mathcal G} }
\def\EE{{\mathcal E}}
\def\KK{{\mathcal K}}
\def\PP{{\mathcal P}}
\def\PPP{{\mathscr P}}
\def\E{{\mathscr E}}
\def\PN{{\mathcal P}}
\def\PPN{{\mathscr P}}
\def\QQ{{\mathcal Q}}
\def\J{J}
\def\Np{{\hat{N}}}
\def\Lp{{\hat{L}}}
\def\Jp{{\hat{J}}}
\def\Pp{{\hat{P}}}
\def\Pip{{\hat{\Pi}}}
\def\Vp{{\hat{V}}}
\def\Ep{{\hat{E}}}
\def\Gp{{\hat{G}}}
\def\Kp{{\hat{K}}}
\def\Ip{{\hat{I}}}
\def\Tp{{\hat{T}}}
\def\Mp{{\hat{M}}}
\def\La{\Lambda}
\def\Ga{\Gamma}
\def\Si{\Sigma}
\def\Upsi{\Upsilon}
\def\Gam{\Gamma}
\def\Gag{{\check{\Gamma}}}
\def\Lap{{\hat{\Lambda}}}
\def\Upsig{{\check{\Upsilon}}}
\def\Kg{{\check{K}}}
\def\ellp{{\hat{\ell}}}
\def\j{j}
\def\jp{{\hat{j}}}
\def\BB{{\mathcal B}}
\def\LL{{\mathcal L}}
\def\MM{{\mathcal U}}
\def\SS{{\mathcal S}}
\def\DD{D}
\def\Dd{{\mathcal D}}
\def\VV{{\mathcal V}}
\def\WW{{\mathcal W}}
\def\OO{{\mathcal O}}
\def\RR{{\mathcal R}}
\def\TT{{\mathcal T}}
\def\AA{{\mathcal A}}
\def\CC{{\mathcal C}}
\def\JJ{{\mathcal J}}
\def\NN{{\mathcal N}}
\def\HH{{\mathcal H}}
\def\XX{{\mathcal X}}
\def\XXX{{\mathscr X}}
\def\YY{{\mathcal Y}}
\def\ZZ{{\mathcal Z}}
\def\CC{{\mathcal C}}
\def\cir{{\scriptscriptstyle \circ}}
\def\circa{\thickapprox}
\def\vain{\rightarrow}
\def\salto{\vskip 0.2truecm \noindent}
\def\spazio{\vskip 0.5truecm \noindent}
\def\vs1{\vskip 1cm \noindent}
\def\fine{\hfill $\square$ \vskip 0.2cm \noindent}
\def\ffine{\hfill $\lozenge$ \vskip 0.2cm \noindent}
\newcommand{\rref}[1]{(\ref{#1})}
\def\beq{\begin{equation}}
\def\feq{\end{equation}}
\def\beqq{\begin{eqnarray}}
\def\feqq{\end{eqnarray}}
\def\barray{\begin{array}}
\def\farray{\end{array}}
\makeatletter \@addtoreset{equation}{section}
\renewcommand{\theequation}{\thesection.\arabic{equation}}
\makeatother
\begin{titlepage}
{~} \vspace{-2.7cm}
\begin{center}
{\huge On power series solutions for the Euler equation,
and the Behr-Ne$\check{\mbox{c}}$as-Wu initial datum}
\end{center}
\vspace{0.5truecm}
\begin{center}
{\large
Carlo Morosi$\,{}^a$, Mario Pernici $\,{}^b$, Livio Pizzocchero$\,{}^c$
({\footnote{Corresponding author}})
} \\
\vspace{0.3truecm} ${}^a$ Dipartimento di Matematica, Politecnico
di Milano,
\\ P.za L. da Vinci 32, I-20133 Milano, Italy \\
e--mail: carlo.morosi@polimi.it \\
\vspace{0.2truecm} ${}^b$ Istituto Nazionale di Fisica Nucleare, Sezione di Milano, \\
Via Celoria 16, I-20133 Milano, Italy \\
e--mail: mario.pernici@mi.infn.it
\\
\vspace{0.2truecm} ${}^c$ Dipartimento di Matematica, Universit\`a di Milano\\
Via C. Saldini 50, I-20133 Milano, Italy\\
and Istituto Nazionale di Fisica Nucleare, Sezione di Milano, Italy \\
e--mail: livio.pizzocchero@unimi.it
\end{center}
\begin{abstract}
We consider the Euler equation
for an incompressible fluid on a three dimensional torus,
and the construction of its solution as a power series
in time. We point out some general facts on this subject,
from convergence issues for the power series
to the role of symmetries of the initial datum. We then turn
the attention to a paper by Behr, Ne$\check{\mbox{c}}$as and Wu
\cite{Nec}; here, the authors chose a very simple
Fourier polynomial as an initial datum for the Euler equation
and analyzed the power series in time for the solution, determining the first 35
terms by computer algebra. Their calculations
suggested for the  series a finite convergence
radius $\tau_3$ in the $H^3$ Sobolev space,
with $0.32 < \tau_3 < 0.35$; they
regarded this as an indication that
the solution of the Euler equation blows up. \par
We have repeated the calculations of \cite{Nec}, using
again computer algebra; the order has been increased
from $35$ to $52$, using the symmetries of the initial
datum to speed up computations. As for $\tau_3$, our results
agree with the original computations of \cite{Nec} (yielding
in fact to conjecture that $0.32 < \tau_3 < 0.33$). Moreover,
our analysis supports the following conclusions: \parn
(a) The finiteness of $\tau_3$ is not at all an indication
of a possible blow-up. \parn
(b) There is a strong indication
that the solution of the Euler equation does not blow up
at a time close to $\tau_3$. In fact,
the solution is likely
to exist, at least, up to a time $\theta_3 > 0.47$. \parn
(c) \Pade analysis gives a rather weak indication
that the solution might blow up at a later time.
\end{abstract}
\vspace{0.0cm} \noindent
\hfill \par
\par \vspace{0.05truecm} \noindent \textbf{AMS 2000 Subject classifications:} 35Q31, 76B03, 35B44, 76M60.
\end{titlepage}
\section{Introduction}
\label{intro}
Let us consider the three-dimensional
Euler equation for a homogeneous incompressible fluid (of unit density) with initial datum $u_0$, i.e.,
$$ {\partial u \over \partial t}  = - u \sc \nabla u - \nabla p~,
\qquad u(x,0) = u_0(x)~.  $$
The unknown is the divergence free
velocity field $(x,t) \mapsto u(x,t)$; we assume periodic boundary conditions, so
$x = (x_1, x_2, x_3)$ ranges in the three dimensional torus
$(\reali/2 \pi \interi)^3$. In the sequel, we often write $u(t)$ for
the function $x \mapsto u(x,t)$. \par
One can try a solution of the above Cauchy problem
in the form of a power series $u(t) = \sum_{j=0}^{+\infty} u_j t^j$
(with $u_j = u_j(x)$); such power series have been the object of rather extensive
investigations.
Morf \textsl{et al} \cite{Fri0}, Frisch \cite{Fri}, Brachet \textsl{et al} \cite{Bra}, Pelz \cite{Pel},
and other authors (see the bibliography of the cited references) have constructed by computer algebra techniques
many terms of the power series for specific initial data, consisting of
simple Fourier polynomials; more precisely, the data analyzed
in these works are the so-called
``Taylor-Green vortex'', and other vortices proposed by Kida \cite{Kida}. The cited authors
have also discussed the possibility of a blow-up (i.e.,
finite-time divergence of $u(t)$) on the grounds of their
computer algebra calculations.
Another initial datum (again a Fourier polynomial)
has been considered by Behr, Ne$\check{\mbox{c}}$as and Wu
\cite{Nec}; these authors have constructed $35$ terms of the power series,
and claimed to have found evidence for a blow-up of the solution; however, in comparison
with the vortices of Taylor-Green and Kida, the Behr-Ne$\check{\mbox{c}}$as-Wu initial
datum has received less attention in the literature. \par
The purpose of the present paper is twofold. \parn
(i) First of all, we wish to point out a number of general facts on the solutions
of the Euler equation and, in particular, on the convergence
of the power series $\sum_{j=0}^{+\infty} u_j t^j$; this is the subject
of Sections \ref{cauc} and \ref{power}. Here we report some
results extracted from the existing literature on the Euler equation
in spaces of analytic functions and/or in Sobolev spaces;
in addition to these results, we present some remarks of ours
and propose a general treatment to discuss the symmetries of the initial
datum and their effects on the solution
of the Euler equation.
We think it is not useless to collect all these theoretical statements
in a unifying framework, suitable for direct application to computer algebra calculations.
\parn
(ii) Our second aim is to reanalyze
the power series for the Behr-Ne$\check{\mbox{c}}$as-Wu
initial datum, both from the theoretical and from the computational
viewpoint; this is the subject of Sections \ref{poneca}, \ref{ourap} and \ref{secpade}.
First of all we apply to the Behr-Ne$\check{\mbox{c}}$as-Wu case
our general setting for the symmetries of the initial datum.
We calculate the symmetry subgroup of the Behr-Ne$\check{\mbox{c}}$as-Wu
datum (that we recognize to be the dihedral group of order $6$;
this group also determines what we call the
pseudo-symmetry space of the datum). \par
With these premises, we present a novel computation of the power series
for the Behr-Ne$\check{\mbox{c}}$as-Wu datum,
based on a Python program written for this purpose; this computation
attains the order $52$.
The Python program uses an exact representation of rational numbers as ratios
of integer, so as not to introduce rounding errors; furthermore,
it employs the symmetries of the initial datum to reduce the
amount of calculations. \par
The results of such computations can be analyzed
using the theoretical framework of Sections \ref{cauc} and \ref{power}. Our conclusions
are the following: \parn
(a) We agree with the estimates of \cite{Nec}, according to which the
power series under consideration has a convergence radius $0.32 < \tau_3 < 0.35$
in the Sobolev space $\Hm{3}$; in fact, our computations suggest $0.32 < \tau_3 < 0.33$.
However, we disagree from the authors of
\cite{Nec} when they interpret the finiteness of $\tau_3$ as indicating
a blow-up of the solution. \parn
(b) On the contrary, we give evidence that the solution
$u(t)$ of the Euler equation exists for $t$ sensibly larger
than $\tau_3$. In fact, analyzing the power series for
the squared Sobolev norm $\| u(t) \|^2_3$, we find a strong indication
for a convergence radius $\theta_3$ such that $0.47 < \theta_3 < 0.50$.
By a general criterion \textsl{\`a la}
Beale-Kato-Majda, this implies that the solution of the Euler equation
exists, at least, up to time $\theta_3$. \par
The final part of our analysis concerns an
alternative approach to estimate $\theta_3$, and
the possibility that $u(t)$ blows up at times larger than $\theta_3$. In connection
with this problem we use the idea (employed in \cite{Bra} \cite{Fri} \cite{Fri0}
\cite{Pel} for different initial data) to construct the
\Pade approximants for the (squared) Sobolev norms
and analyze their singularities.
In particular, we construct the diagonal \Pade approximants
$[p/p](t)$ for $\| u(t) \|^2_3$, up to $p=26$. For
most of them the complex singularities of minimum modulus
have modulus $\simeq 0.5$; this fact yields new
evidence for the previous estimate on $\theta_3$.
Moreover, most of these \Pade approximants
have real singularities, distributed rather erratically; analyzing
them in terms of mean value and variance, we obtain a somehow weak indication
that: \parn
(c) $u(t)$ might blow up for $t \vain T^{-}$ (and $t \vain (-T)^{+}$),
for some $T$ such that $0.56 < T < 0.73$. \par
The blow-up problem can be studied as well in terms
of D-log \Pade approximants; these do not give a clear
indication supporting conjecture (c), as briefly explained
at the end of the paper.
In general, much caution is recommended
about the Euler equation and
blow-up predictions via \Pade analysis: for example, in the case
of the Taylor-Green vortex the \Pade approximants exhibit
real singularities \cite{Fri} \cite{Fri0}, but the numerical solution of the Euler
equation by spectral methods raises doubts on
the actual existence of a blow-up \cite{Bra1} \cite{Bra2}.
\parn
\textbf{Connections with other works.}
Concluding this Introduction, to put the subject of this
paper into a wider perspective we wish to mention that there are
general methods of functional analysis to obtain
quantitive lower bounds on the time of existence $T$
of the solution of the Euler (or Navier-Stokes)
Cauchy problem, from the \textsl{a posteriori}
analysis of an approximate solution; such lower
bounds are certain (i.e., non conjectural). \par
Derivations of such \textsl{a posteriori} lower bounds
have been given in \cite{Che} \cite{appro} \cite{appeul}. The last
of these works gives an algorithm to obtain
these lower bounds analyzing any approximate solution
of the Euler (or Navier-Stokes)
Cauchy problem via a suitable differential inequality,
called therein the ''control inequality''.
\par
Again in \cite{appeul}, a preliminary analysis of the Euler (and Navier-Stokes)
equations with the Behr-Ne$\check{\mbox{c}}$as-Wu initial
datum has been performed, using for the solution
a Galerkin approximation with very few Fourier modes. This
approximant, combined with the control inequality,
gives for the Euler equation with this datum a (poor, but
certain) lower bound $T > 0.066$ for the time of existence
in $H^3$ (the same approach, applied to the Navier-Stokes
equations, grants $T = + \infty$ when the viscosity
coefficient is above an explicit threshold).
We plan to continue in future works the analysis of the
Behr-Ne$\check{\mbox{c}}$as-Wu intial datum, combining
the control equation of \cite{appeul} with approximation methods
based on extensive automatic computations such as the ones
presented in this paper.
\section{The Cauchy problem for the Euler equation on a torus}
\label{cauc}
\textbf{Preliminaries.}
If $a = (a_s), b = (b_s)$ are elements of $\reali^3$ or $\complessi^3$,
we intend $a \sc b := \sum_{s=1}^3 a_s b_s$.
We indicate with $\overline{\phantom{b}}$ the complex conjugate (and
we let it act componentwise on elements of $\complessi^3$); we put
$|a| := \sqrt{\overline{a} \sc a} = \sqrt{\sum_{s=1}^3 |a_s|^2}$.
\par
The Cauchy problem for the incompressible Euler equation is
\beq {\partial u \over \partial t}  = - u \sc \nabla u - \nabla p~,
\qquad u(x,0) = u_0(x)~, \label{eulp} \feq
where: $u= u(x, t)$ is the divergence free velocity field;
the space variables $x = (x_s)_{s=1,2,3}$ belong to the torus
$\Tt := (\reali/2 \pi \interi)^3$;
$(u \sc \nabla u)_r := \sum_{s=1}^3 u_s \partial_s u_r$ ($r=1, 2, 3$);
$p = p(x,t)$ is the pressure; $u_0 = u_0(x)$ is the initial
datum. As well known, the pressure can be
eliminated from \rref{eulp} using
the Leray projection $\LP$ onto the space of divergence free vector fields;
this allows to rewrite the evolution equation in \rref{eulp} as
${\partial u/\partial t}  = - \LP(u \sc \nabla u)$. In this way, we obtain
for the Cauchy problem the final form
\beq {\partial u \over \partial t}  = \PPP(u,u)~, \qquad u(\, . \, ,0) = u_0~,\label{eul} \feq
where we have written $\PPP$ for the bilinear map sending two
(sufficiently regular) vector fields $v, w : \Tt \vain \reali^3$ into the vector field
\beq \PPP(v,w) := - \LP(v \sc \nabla w)~. \feq
In this framework, it is convenient to associate to
a vector field $v : \Tt \vain \reali^ 3$ the Fourier
components $v_k := (2 \pi)^{-3} \int_{\Tt} d x \, e^{-i k \sc x} v(x) \in \complessi^3$, so that
\beq v(x) = \sum_{k \in \Zt} v_k e^{i k \sc x}~. \feq
Due to the reality of $v$, we have $v_{-k} = \overline{v_k}$, and $v$ is
divergence free iff $k \sc v_k=0$ for all $k$. With $\PPP$ as above
and $v, w$ two vector fields, the Fourier
components of $\PPP(v, w)$ are
\beq \PPP(v, w)_k = - i
\sum_{h \in \Zt} v_h \sc (k-h) \, \LP_k w_{k-h} \label{repp} \feq
where $\LP_k : \complessi^ 3 \vain \complessi^3$ is the
projection on the orthogonal complement of $k$ ($\LP_k c :=
c - (k \sc c) k/|k|^2$ if $k\neq 0$; $\LP_0 c := c$). \par
In the above, we have introduced the setting for the Euler equation
in an informal way; to go on, it is necessary to
specify the functional spaces to which the velocity fields
(at any time) are supposed to belong.
\par
The expression ``a vector field $\Tt \vain \reali^3$'' can be understood,
with very wide generality, as
``an $\reali^3$-valued distribution on $\Tt$'' (see, e.g.,
\cite{accau}); we write
$D'(\Tt, \reali^3) \equiv D'$ for the space of such distributions.
Any $v \in D'(\Tt, \reali^3)$ can be differentiated in the
distributional sense and has a (weakly convergent) Fourier expansion with coefficients
$v_k \in \complessi^3$, such that $\overline{v_k} = v_{-k}$. \par
To construct the full setting for the Euler equation, one must
confine the attention to much smaller functional spaces of
vector fields. For our purposes, two cases
are important: \parn
(i) The Sobolev space $\Hm{n}$ of zero mean, divergence free vector
fields of any order $n \in [0,+\infty)$. This is defined in terms
of the space $L^2(\Tt, \reali^3) \equiv L^2$ of square integrable
vector fields $v: \Tt \vain \reali^3$, equipped with the inner
product $\la v | w \ra_{L^2} := (2 \pi)^{-3} \int_{\Tt} v(x) \sc w(x) d x$
and with the induced norm $\| v \|_{L^2} = (2 \pi)^{-3/2} \sqrt{\int_{\Tt} |v(x)|^2 d x}$
(note the term $(2 \pi)^{-3}$ in the inner product, used systematically in the sequel).
By definition,
\parn
\vbox{
\beq \Hm{n}(\Tt, \reali^3) \equiv \Hm{n} := \Big\{ v \in D'~|~
\sqrt{-\Delta}^{\,n} v \in L^2~,
\int_{\Tt} \! \! v~ d x = 0, ~\dive \,v = 0 \Big\}  \feq
$$ = \Big\{ v \in D'~|~
\sum_{k \in \Zt} |k|^{2 n} |v_k|^2 < + \infty, v_0 = 0,
k \sc v_k = 0  \Big\}~. $$
}
\noindent
(In the above $\sqrt{-\Delta}^n$
indicates the power of order $n/2$ of minus the Laplacian;
by definition $(\sqrt{-\Delta}^n v)_k = |k|^n v_k$
for each $v \in D'$. Note that $H^n \subset L^2$
for all $n \geqs 0$.)
$H^n$ is a Hilbert space with the inner product
\beq \la v | w \ra_{n} := \la \sqrt{-\Delta}^n v | \sqrt{-\Delta}^n w \ra_{L^2} = \
\sum_{k \in \Zt} |k|^{2 n} \overline{v_k} \sc w_k~, \feq
inducing the norm
\beq \| v \|_n = \| \sqrt{-\Delta}^{\,n} v \|_{L^2} =
\sqrt{\sum_{k \in \Zt} |k|^{2 n} |v_k|^2}~. \feq
It is known that $\PPP$ sends continuously $H^n \times H^{n+1}$ into $H^n$,
for all $n \in (3/2,+\infty)$. \parn
(ii) The space of $C^\om$ (i.e., analytic) zero mean, divergence free vector fields
on $\Tt$; this is
\beq \Aa(\Tt, \reali^3) \equiv \Aa := \Big\{ v \in C^\om(\Tt, \reali^3)~|~
\int_{\Tt} \! \! v~ d x = 0, ~\dive v = 0~\Big\}  \label{29} \feq
$$ = \Big\{ v \in D'~|~\liminf_{k \in \Zt, \, k \vain \infty}
|v_k|^{-{1 \over |k_1| + |k_2| + |k_3|}} > 1, v_0 = 0,~
k \sc v_k = 0 \Big \}~ $$
(intending $0^{-{1 \over |k_1| + |k_2| + |k_3|}} := +\infty$.
The Fourier representation in \rref{29}
mimics the description of analytic functions on the torus in \cite{Mori}, which is also
a useful reference for what follows).
One has \parn
\vbox{
\beq \Aa = \cup_{\rho \in (1,+\infty)} \Aa_\rho~, \label{defaa} \feq
$$ \Aa_\rho := \Big\{ v \in D'~|~\liminf_{k \in \Zt, \, k \vain \infty}
|v_k|^{-{1 \over |k_1| + |k_2| + |k_3|}} > \rho,~ v_0 = 0,~
k \sc v_k = 0 \Big\}~; $$
}
\noindent
each $\Aa_\rho$ is a vector subspace of $\Aa$.
Let us introduce the annulus $K_\rho := \{ z \in \complessi~|~1/\rho \leqs |z | \leqs \rho\}$ and
its power $K^3_\rho := \{ z = (z_1, z_2, z_3) \in \complessi^3~|~z_1, z_2, z_3 \in K_\rho \}$.
For $v \in \Aa_\rho$, the series $\sum_{k \in \Zt} v_k z^k$ converges in $\complessi^3$ for each
$z \in K^3_\rho$ (intending $z^k := z_1^{k_1} z_2^{k_2} z_3^{k_3})$;
the function $z \mapsto \sum_{k \in \Zt} v_k z^k$ is holomorphic on the inner part
of $K^3_\rho$ and continuous on $K^3_{\rho}$, so we can define
\beq \bbb v \bbb_\rho := \sup_{z \in K^3_\rho} \left|\sum_{k \in \Zt} v_k z^k \right|~. \feq
$\bbb~\bbb_\rho$ is a norm on $\Aa_\rho$ and makes it a Banach space. One equips $\Aa$ with
the inductive limit topology of the collection of Banach spaces
$\{(\Aa_\rho, \bbb~\bbb_\rho)~|~\rho \in (1,+\infty) \}$:
this is the finest locally convex topology on $\Aa$ making continuous each
embedding $\Aa_\rho \hookrightarrow \Aa$. (Besides \cite{Mori},
see \cite{Trev} for the general theory of inductive limits.)
$\Aa$ is continuously embedded into each Sobolev space $\Hm{n}$;
the map $\PPP$ is continuous from $\Aa \times \Aa$ to $\Aa$.
\salto
\textbf{Basic results on local existence and uniqueness.}
We start from the Sobolev framework, choosing
\beq n \in (5/2, + \infty)~. \feq
In the sequel, an $\Hm{n}$-solution of the Euler equation, or of the Euler Cauchy
problem, means a map
\beq u \in C((-\Ti,T),\Hm{n}) \cap C^1((-\Ti,T),\Hm{n-1}) \label{maxi1} \feq
($\Ti, T \in (0,+\infty]$) fulfilling the Euler equation, or its Cauchy
problem with a suitable initial condition $u(0) = u_0$.
The following statement is well known:
\begin{prop}
\label{prop1}
\textbf{Proposition.} For $n$ as above and any initial datum $u_0
\in H^n$, the following holds.
\parn
(i) The Cauchy problem \rref{eul} has a unique maximal (i.e., not extendable) $\Hm{n}$-solution
$u$ of domain $(-\Ti, T)$,
for suitable $T = T(u_0), \Ti = \Ti(u_0) \in (0,+\infty]$. \parn
(ii) (Beale-Kato-Majda criterion, Sobolev version). If $T < + \infty$, one has
\beq \int_{0}^T \! d t \, \| u(t) \|_n = + \infty~, \label{blowhni} \feq
a fact implying
\beq \limsup_{t \vain T^{-}} \| u(t) \|_n  = + \infty~. \label{blowhn} \feq
Similar results hold if $\Ti < + \infty$, considering the integral from $-\Ti$ to $0$ and
the limit for $t \vain (- \Ti)^{+}$.
\end{prop}
\textbf{Proof.} (i) See \cite{BKM} \cite{Kat2}. \parn
(ii) See \cite{BKM}. Indeed, here it is shown that $T < + \infty$ implies
$\int_{0}^T d t \| \rot u(t) \|_{L^\infty} = +\infty$; however,
$\| \rot u(t) \|_{L^{\infty}} \leqs \mbox{const.} \| u(t) \|_n$
by the Sobolev imbedding inequalities, whence Eq. \rref{blowhni} and its obvious
consequence \rref{blowhn}. The behavior
of $u$ at time $-\Ti$ is analyzed similarly. \fine
\par
If $T < +\infty$, the solution $u$ is said to blow up at time $T$.
Similarly, if $\Ti < + \infty$
we say that $u$ blows up at $-\Ti$. Many statements presented in the sequel
on the possibility of blow-up at $T$ have obvious reformulations regarding $-\Ti$.
\begin{rema}
\label{remabkm}
\textbf{Remark.} The Beale-Kato-Majda criterion \rref{blowhni}
yields the following statement, in case of blow-up with a power law:
\beq \mbox{if}~ \| u(t) \|_n \sim {U \over (T - t)^{\alpha}}
\quad \mbox{for $t \vain T^{-}$ (with $U, \alpha > 0$)}, \quad \mbox{then $\alpha \geqs 1$}~. \label{teorkato} \feq
In the case of the Euler equation on $\reali^3$, it was recently
shown in \cite{Chen}, Theorem 1.3 that the blow-up at $T$ implies the following,
for any $n > 5/2$:
\beq \| u(t) \|_n \geqs {U \over (T - t)^{2 n/5}}
\quad \mbox{for $t$ close to $T$}~,
\qquad U = U_n(\| u_0 \|_{L^2})~.
\label{teorest} \feq
This estimate might hold as well for the framework
of the present paper, i.e., for the Euler equation
on the torus $\Tt$ (however, the extendability of \rref{teorest}
to $\Tt$ is immaterial for the purposes of this paper).  \fine
\end{rema}
\par
Let us pass to the $C^\om$ ($=$ analytic) framework; what follows
assumes some general notions from the theory
of analytic functions from $\reali$ to locally convex spaces,
for which we refer to \cite{Bour}  \S 3. Let $\Aa$ be the space \rref{29};
in the sequel, an $\Aa$-solution of the Euler equation, or of the Euler Cauchy
problem, means a map
\beq u \in C^\om((-\Ti,T),\Aa)~\label{uanal} \feq
($\Ti, T \in (0,+\infty]$) fulfilling the Euler equation, or its Cauchy
problem with a suitable initial condition $u(0) = u_0$.
Let us report a known result.
\begin{prop}
\label{prop2}
\textbf{Proposition.} For any initial datum $u_0
\in \Aa$, the following holds.
\parn
(i) Problem \rref{eul} has a unique maximal (i.e., non extendable) $\Aa$-solution
of domain $(-\Ti, T)$,
for suitable $T = T(u_0), \Ti = \Ti(u_0) \in (0,+\infty]$. \parn
(ii) For any $n \in (5/2, + \infty)$, this coincides with the maximal
$\Hm{n}$-solution of the Cauchy problem with the same datum (and thus, if $T < + \infty$,
it fulfills Eqs. \rref{blowhni} \rref{blowhn}; a similar result holds if $\Ti < + \infty$).
\end{prop}
\textbf{Proof.} (i) See \cite{Bau}, Theorem III.2, page 264 (this is a result of
existence and uniqueness on sufficiently small time intervals,
from which one infers via standard arguments existence and uniqueness
of the maximal solution). \parn
(ii) See \cite{Bar}, especially Remark 2.1, page 414.
\fine
Assuming again $u_0 \in \Aa$, and choosing any $n \in [0,+\infty)$,
we conclude with two remarks. \parn
(i) By the continuous embedding of
$\Aa$ into $\Hm{n}$, the function $u$ of the last proposition
is also in $C^\om((-\Ti,T), \Hm{n})$. \parn
(ii) Consider the function
\beq (-\Ti,T) \vain \reali~, \qquad t \mapsto \| u(t) \|_n^2~. \label{nun} \feq
This is in $C^\om((-\Ti,T), \reali)$,
being the composition of the analytic
function $u : (-\Ti,T)$ $\vain \Hm{n}$ with the continuous quadratic
function $\| ~\|^2_n : \Hm{n} \vain \reali$.
\salto
\salto
\textbf{Symmetries of the Euler equation.} Let us consider the \textsl{octahedral group}
$O_h$, formed by the orthogonal $3 \times 3$ matrices with integer entries:
\beq O_h := \{ S \in Mat(3 \times 3, \interi)~|~S^\t S = \uno_3 \}~. \feq
In fact, the entries of any such matrix have $-1,0$ and $1$ as the only possible values;
furthermore, a $3 \times 3$ matrix $S$ belongs to $O_h$ if and only if
\parn
\vbox{
\beq S = \mbox{diag}(\ep_1,\ep_2,\ep_3)  Q(\sigma) \feq
$$ \ep_s \in \{\pm 1\}~(s=1,2,3) \, ;~ Q(\sigma)~
\mbox{the matrix of the permutation $\sigma : \{1,2,3\} \vain \{1,2,3\}$}; $$}
more precisely, $Q(\sigma)$ is the matrix such that $(Q(\sigma) c)_s = c_{\sigma(s)}$
for all $c \in \complessi^3$, $s\in \{1,2,3\}$. There are $2^3 = 8$ possible choices
for the signs $\ep_i$ and $3! = 6$ choices for $\sigma$, so $O_h$ has $8 \times 6 = 48$
elements. Clearly, each $S \in O_h$ sends $\Zt$ into itself. \par
To go on, let us denote with $\ohz$ the Cartesian product
$\ohperz$, viewed as a group with the composition law defined by
({\footnote{This is the semidirect product
of the groups $O_h$ and $\Tt$ with respect to the natural homomorphism
$O_h \vain Aut(\Tt)$ sending $S \in O_h$ into
the map $b \mapsto S b$, an automorphism of $\Tt$.}})
\beq (S,a) (U,b) := (S U, a + S b) \quad (S, U  \in O_h~;~ a, b \in \Tt)~. \label{prcirc} \feq
Of course, the unit of this group is $(\uno,0)$ (with $\uno$ the
identity $3 \times 3$ matrix); the inverse of a pair $(S,a)$ is $(S,a)^{-1} = (S^\t, - S^\t a)$.
To any element $(S,a)$ of $\ohz$ is associated a ``rototranslation''
\beq \E(S,a) : \Tt \vain \Tt~, \qquad x \mapsto \E(S,a)(x) := S x + a~, \feq
and one checks that the mapping $(S,a) \mapsto \E(S,a)$ is a group homomorphism
between $\ohz$ and the group of diffeomorphisms of $\Tt$ into itself
(with the usual composition).  \par
Now, we take a vector field
$v$ in $H^n$ (or in $\Aa$) and an element $(S,a)$ of the group $\ohz$.
We can construct the push-forward $\E_{*}(S,a) v$
of $v$ along the mapping $\E(S,a)$; this is the vector field
in $H^n$ (or in $\Aa$), given by
\beq \E_{*}(S,a) v : \Tt \vain \reali^3~, \qquad x \mapsto
(\E_{*}(S,a) \, v)(x) = S v(S^\t(x - a))~. \label{push} \feq
One easily checks that Eq. \rref{push} actually defines a vector field in $H^n$ (or in $\Aa$), with
Fourier components
\beq (\E_{*}(S,a) \, v)_k = e^{-i a \sc k} S v_{\scriptscriptstyle{ S^\t k}} \qquad (k \in \Zt)~. \label{pushf} \feq
Let us write $\E_{*}(S,a)$ for the map $v \in H^n \mapsto \E_{*}(S,a) v$;
this is a linear map of $H^n$ into itself,
preserving the inner product $\la~|~\ra_n$, so it is
in the group $O(H^n)$ of orthogonal operators of the Hilbert space $H^n$ into itself.
The mapping
\beq \E_{*}: \ohz \vain O(H^n)~, \qquad (S, a) \mapsto \E_{*}(S,a) \feq
is a injective group homomorphism, i.e., a faithful orthogonal representation of the group $\ohz$ on
the real Hilbert space $H^n$.
Alternatively, let us write
$\E_{*}(S,a)$ for the map $v \in \Aa \mapsto \E_{*}(S,a) v$;
this is in the space $\mbox{Iso}(\Aa)$ of linear and topological isomorphisms of $\Aa$ into itself.
The map
\beq \E_{*}: \ohz \vain \mbox{Iso}(\Aa)~, \qquad (S, a) \mapsto \E_{*}(S,a) \feq
is an injective group homomorphism, i.e., a faithful linear representation of the group $\ohz$ on
the topological vector space $\Aa$. \par
Let us relate the previous constructions to the bilinear map $\PPP$ of the
Euler equation. From the Fourier representations \rref{repp} \rref{pushf}, one easily infers
\beq \PPP(\E_{*}(S,a) \, v, \E_{*}(S,a) \, w) = \E_{*}(S,a)\,  \PPP(v,w) \label{pinv} \feq
for all $v \in H^n$, $w \in H^{n+1}$ with $n > 3/2$ (and, in particular,
for all $v, w \in \Aa$).
Let us outline the implications of \rref{pinv} about the solutions of
the Euler equation. In the rest of the paragraph,
the term ``solution'' either means an $H^n$-solution $(n> 5/2)$
or an $\Aa$-solution, and the initial datum $u_0$ is chosen consistently in $H^n$ or in $\Aa$.
From \rref{pinv} one infers the following,
for each $(S,a) \in \ohz$: \parn
(i) If $u : t \in (-\Ti, T) \mapsto u(t)$ is a solution of the Euler equation,
we have two more solutions
\beq \E_{*}(S,a) u: ~t \in (-\Ti, T) \mapsto \E_{*}(S,a) u(t), \feq
\beq - \E_{*}(S,a) u(-\cdot): ~t \in (-T,\Ti) \mapsto - \E_{*}(S,a) u(-t). \feq
(ii) If $u : t \in (-\Ti, T) \mapsto u(t)$ is the maximal solution of the
Euler Cauchy problem with datum $u_0$, then $\E_{*}(S,a) u$
is the maximal solution with datum $\E_{*}(S,a) u_0$ and
$-\E_{*}(S,a) u(-\cdot)$ is the maximal solution with datum
$-\E_{*}(S,a) u_0$. \parn
(iii) Let us denote again with $u : t \in (-\Ti, T) \mapsto u(t)$ the maximal solution of
the Cauchy problem with datum $u_0$. Then,
\beq \E_{*}(S,a) u_0  = u_0~\Rightarrow~
\E_{*}(S,a) u(t) = u(t) \qquad \mbox{for $t \in (-\Ti,T)$}. \label{rs2} \feq
\beq - \E_{*}(S,a) u_0  = u_0~\Rightarrow~\Ti=T,~
-\E_{*}(S,a) u(-t) = u(t) \qquad \mbox{for $t \in (-T,T)$}. \label{rs22} \feq
The verification of statements (i)(ii) is straightforward. After this,
the implication \rref{rs2} in (iii) follows
noting that $\E_{*}(S,a) u$ and $u$
are maximal solutions of the Cauchy problem with the same datum
$\E_{*}(S,a) u_0 = u_0$. Similarly,
the implication \rref{rs22} follows noting that
$-\E_{*}(S,a) u(-\cdot)$ and $u$
are maximal solutions of the Cauchy problem with the same datum
$-\E_{*}(S,a) u_0 = u_0$. \par
Considering the maximal solution $u$ for a datum $u_0$ in $H^n$ ($n > 5/2$), and
recalling that any transformation $\E_{*}(S,a)$ preserves
the $H^n$ norm, we also obtain from \rref{rs22} the following:
\beq - \E_{*}(S,a) u_0  = u_0~\Rightarrow~\Ti=T,~
\| u(-t) \|_n = \| u(t) \|_n \qquad \mbox{for $t \in (-T,T)$}. \label{rs22inv} \feq
The results in (iii) suggest to consider, for a given datum $u_0$ in $H^n$ or $\Aa$,
the \textsl{symmetry subgroup}
\beq \Isot(u_0) := \{ (S,a) \in \ohz~|~\E_{*}(S,a) u_0 = u_0 \} \label{isot} \feq
and the \textsl{pseudo-symmetry space}
\beq \PIsot(u_0) := \{ (S,a) \in \ohz~|~-\E_{*}(S,a) u_0 = u_0 \} \label{pisot} \feq
(the first one, being a subgroup of $\ohz$, contains at least the identity
element $(\uno, 0)$; the second one might be the empty set. The term
''isotropy group'', often employed in place of ''symmetry group'', will
not be used in this paper). \parn
Let us consider the maximal solution $u$ of the Cauchy problem with a datum
$u_0$ (contained in $H^n$ for some $n > 5/2$);
from Eqs. \rref{rs2}-\rref{rs22inv}, we readily obtain the
following:
\beq \E_{*}(S,a) u(t) = u(t) \qquad \mbox{for all $(S,a) \in \Isot(u_0),
t \in (-\Ti,T)$}~; \feq
\beq \PIsot(u_0) \neq \emptyset\quad \Rightarrow \quad \Ti=T,~
-\E_{*}(S,a) u(-t) = u(t), ~\| u(-t) \|_n = \| u(t) \|_n\feq
$$\quad\mbox{for $(S,a) \in \PIsot(u_0),
t \in (-T,T)$}~. $$
For future use, let us introduce the \textsl{reduced symmetry subgroup} and the
\textsl{reduced pseudo-symmetry space} of the
datum $u_0$, which are
\beq \Isott(u_0) := \{ S \in O_h~|~\E_{*}(S,a) u_0 = u_0
~\mbox{for some $a \in \Tt$} \}~, \label{isott} \feq
\beq \PIsott(u_0) := \{ S \in O_h~|~-\E_{*}(S,a) u_0 = u_0
~\mbox{for some $a \in \Tt$} \}~. \label{pisott} \feq
Let us observe that the set theoretical unions
$\Isot(u_0) \cup \PIsot(u_0)$ and $\Isott(u_0) \cup \PIsott(u_0)$
are subgroups of $\ohz$ and $O_h$, respectively. \par
As a final remark, useful for the sequel, let us
consider the pair $(-\uno, 0) \in \ohz$, noting that $\EE(-\uno,0)$
is the space reflection: $\EE(-\uno,0)(x) = - x$ for all $x \in \Tt$. One easily checks
that
\beq (-\uno, 0) \in \PIsot(u_0)
~\Leftrightarrow~ \PIsot(u_0)
 = \Isot(u_0) (-\uno,0) = \{ (-S,a)~|~(S, a) \in \Isot(u_0) \}
\label{pseudorifl} \feq
(where
$\Isot(u_0) (-\uno,0)$ stands for the set
$\{ (S,a)(-\uno,0) ~|~ (S, a) \in \Isot(u_0) \}$;
the last equality rests on the identity
$(S,a) (-\uno,0) = (-S,a)$).
\section{Power series in time for the Euler Cauchy problem}
\label{power}
Throughout this section, we consider the Euler Cauchy problem with initial datum $u_0 \in \Aa$.
\salto
\textbf{Setting up a power series for the solution.}
Let us try to build the solution of the Euler Cauchy problem as a power series
\beq t \mapsto \sum_{j=0}^{\infty} u_j t^j \label{espa} \feq
with coefficients $u_j \in \Aa$, whose convergence has to be discussed
later. The zero order term in this expansion
is the initial datum $u_0$; to determine the other coefficients
$u_\j \in \Aa$, it suffices to substitute the expansion \rref{espa} into
the Euler equation \rref{eul}, and to require equality of the coefficients
of the same powers of $t$ in both sides: in this way, one easily
obtains the recurrence relation
\beq u_j = {1 \over j} \sum_{\ell=0}^{j-1} \PPP(u_\ell, u_{j - \ell - 1})
\qquad (j=1,2,3,...)~. \label{recur} \feq
When applying this recurrence relation for the $u_j$'s
it can be useful to represent the bilinear map $\PPP$ in terms
of Fourier coefficients, as in Eq. \rref{repp}. This is especially
useful if the initial datum $u_0$ is a Fourier polynomial, i.e.,
if $u_{0 k} \neq 0$ only for finitely many modes $k$. In this
case, all the iterates $u_j$ ($j=1,2,3,..$) are as well
Fourier polynomials, and the implementation of \rref{recur}
via the Fourier representation \rref{repp} always involves sums
over finitely many modes. \par
In the next section, a large part of our attention will be devoted
(for a specific datum $u_0$) to the partial sums
\beq u^{(N)}(t) := \sum_{j=0}^N u_j t^j~, \feq
($N=0,1,2,...$) and to the (squared) Sobolev norms
\beq \| u^{(N)}(t) \|^2_n  = \sum_{k \in \Zt} |k|^{2 n} |u^{(N)}_k(t)|^2~. \feq
\salto
\textbf{Symmetry considerations.} Let us consider the symmetry subgroup $\Isot(u_0)$
or the pseudo-symmetry space $\PIsot(u_0)$, see Eqs. \rref{isot} \rref{pisot}.
Using the recursive definition \rref{recur} of $u_j$ with the
invariance property \rref{pinv} of $\PPP$, one easily checks the following,
for any $j \in \{0,1,2,...\}$:
\beq \E_{*}(S,a) u_j = u_j \qquad \mbox{for all $(S,a) \in \Isot(u_0)$}; \label{equal} \feq
\beq - \E_{*}(S,a) u_j = (-1)^{j} u_j \qquad \mbox{for all $(S,a) \in \PIsot(u_0)$}. \label{pequal} \feq
Of course, the last two equations imply the following, for all $N \in \{0,1,2,...\}$, $t \in \reali$
and $n \in [0,+\infty)$:
\beq \E_{*}(S,a) u_N(t) = u_N(t) \qquad \mbox{for $(S,a) \in \Isot(u_0)$}; \label{nequal} \feq
\beq - \E_{*}(S,a) u_N(t) = u_N(-t) \qquad \mbox{for $(S,a) \in \PIsot(u_0)$}; \label{npequal} \feq
\beq \| u_N(t) \|_n = \| u_N(-t) \|_n \qquad \mbox{if $\PIsot(u_0) \neq \emptyset$} \label{norequal} \feq
(Eq. \rref{norequal} is a consequence of Eq. \rref{npequal} and of the invariance
of $\|~\|_n$ under the transformation $\E_{*}(S,a)$). \par
Due to the Fourier representation \rref{pushf} for $\E_{*}(S,a)$,
the equality \rref{equal} reads
$e^{- i a \sc k}$ $ S u_{j, S^\t k}$ $= u_{j, k}$ or, equivalently,
\beq u_{j, S k} = e^{- i a \sc S k} S u_{j, k} \qquad \mbox{for $k \in \Zt$, $(S,a) \in \Isot(u_0)$}~;
\label{equal2} \feq
similarly, Eq. \rref{pequal} is equivalent to the statement
\beq u_{j, S k} = (-1)^{j+1}
e^{- i a \sc S k} S u_{j, k} \qquad \mbox{for $k \in \Zt$, $(S,a) \in \PIsot(u_0)$}~.
\label{pequal2} \feq
In typical applications of the recursion scheme \rref{recur}, where
$u_0$ is a Fourier polynomial as well as its iterates $u_j$, Eqs. \rref{equal2}
\rref{pequal2} can be used to speed up the computation of the Fourier components
of the $u_{j}$'s; in fact, at any given order $j$, after computing
a Fourier component $u_{j, k}$ we immediately obtain from
the cited equations the components $u_{j, S k}$ for all $S$ in the reduced
subgroup or subspace $\Isott(u_0)$, $\PIsott(u_0)$.
\salto
\textbf{Convergence of the power series in $\boma{\Aa}$.}
From now on, we intend
\beq \tau := \mbox{convergence radius of the series $\sum_{j=0}^{\infty} u_j t^j$ in $\Aa$}~. \feq
Furthermore,
\beq u : t \in (-\Ti,T) \mapsto u(t)~ \mbox{is the maximal $\Aa$-solution of
the Cauchy problem} \feq
(recall that, for any $n > 5/2$, $u$ is also the maximal $\Hm{n}$-solution). We note that
\beq 0 < \tau \leqs \Ti \wedge T~, \qquad u(t) = \sum_{j=0}^{\infty} u_j t^j
\quad \mbox{in $\Aa$, for $t \in (-\tau,\tau)$} \feq
(with $\wedge$ indicating the minimum).
In fact: being analytic, $u$ admits a power series representation
in a neighborhood of zero; this necessarily coincides with the
series \rref{espa}, whose convergence radius $\tau$ is thus nonzero
and fulfills $(-\tau,\tau) \subset (-\Ti,T)$.
\salto
\textbf{Convergence of the power series in $\boma{\Hm{n}}$.}
After fixing $n \in [0,+\infty)$, let us discuss the series
\rref{espa} in the Sobolev space $\Hm{n}$. To this purpose, we put
\beq \tau_n := \mbox{convergence radius of the series
$\sum_{j=0}^{\infty} u_j t^j$ in $\Hm{n}$}~; \label{taun} \feq
the root test gives
\beq \tau_{n} = \liminf_{j \vain + \infty} \| u_j \|_n^{-1/j} \label{roote} \feq
(intending $0^{-1/j} := +\infty$).
With $\tau, \Ti, T, u$ as before, we claim that
\beq \tau \leqs \tau_n~\mbox{and}~ u(t) = \sum_{j=0}^{\infty} u_j t^j
~\mbox{in $\Hm{n}$, for $t \in (-\Ti \wedge \tau_n, T \wedge \tau_n)$}~\feq
(where $- \Ti \wedge \tau_n$ is
the opposite of the minimum $\Ti \wedge \tau_n$).
In fact: the series $\sum_{j=0}^{\infty} u_j t^j$ converges
to $u(t)$ in $\Aa$, for $t \in (-\tau, \tau)$; by the continuous embedding
$\Aa \hookrightarrow \Hm{n}$ this series converges to $u(t)$
in $\Hm{n}$ as well, at least for $t \in (-\tau, \tau)$;
thus $\tau \leqs \tau_n$. Moreover the functions $u : (-\Ti,T)
\vain \Hm{n}$ and $t \in (-\tau_n, \tau_n) \mapsto \sum_{j=0}^{\infty} u_j t^j \in H^n$
are analytic and coincide on $(-\tau,\tau)$; so, by the analytic continuation
principle, these functions coincide on the intersection of their domains
which is $(-\Ti \wedge \tau_n, T \wedge \tau_n)$.
Let as add a stronger claim:
\beq \mbox{if $n > \dd{5 \over 2}$}, \qquad
\tau \leqs \tau_n \leqs \Ti \wedge T~\mbox{and}~ u(t) = \sum_{j=0}^{\infty} u_j t^j
~\mbox{in $\Hm{n}$, for $t \in (-\tau_n,\tau_n)$}~. \label{claimm} \feq
In fact, the function $t \in (-\tau_n, \tau_n) \mapsto \sum_{j=0}^{\infty} u_j t^j$\
is in  $C((-\tau_n,\tau_n),\Hm{n})$ $\cap$ $C^1((-\tau_n,\tau_n),$ $\Hm{n-1})$
and solves the Euler Cauchy problem, so it is a restriction
of the maximal $\Hm{n}$-solution, which is $u$ of domain $(-\Ti,T)$; this gives the relations
$\tau_n \leqs \Ti \wedge T$ and $u(t) = \sum_{j=0}^{\infty} u_j t^j$
in $\Hm{n}$, for $t \in (-\tau_n, \tau_n)$.
\salto
\textbf{Power series for the Sobolev norms of the
solution.} Let us choose $n \in [0,+\infty)$.
The squared norm $\| \sum_{j=0}^{+\infty} u_j t^j \|^2_n =
\la \sum_{j=0}^{+\infty} u_j t^j  | \sum_{j=0}^{+\infty} u_j t^j \ra_n $ has the formal expansion
\beq \| \sum_{j=0}^{+\infty} u_j t^j \|^2_n = \sum_{j=0}^{+\infty} \nu_{ n j} t^j,
\qquad \nu_{n j} := \sum_{\ell=0}^j \la u_\ell | u_{j - \ell} \ra_n \in \reali~; \label{formal} \feq
for future use we remark that
({\footnote{Let us propose a proof of \rref{rec3}, based directly on
the definition \rref{formal} of $\nu_{n j}$.
If $\PIsot(u_0)$ has at least one element
$(S,a)$, from \rref{pequal} and from the invariance of $\la~|~\ra_n$ under
any transformation $\E_{*}(S, a)$ we obtain that, for each $\ell \in \{0,...,j\}$,
$\la u_\ell | u_{j - \ell} \ra_n = \la (-1)^\ell \E_{*}(S, a) u_\ell | (-1)^{j-\ell}
\E_{*}(S, a) u_{j - \ell} \ra_n = (-1)^{j} \la u_\ell | u_{j - \ell} \ra_n$,
whence $\nu_{n j} = (-1)^j \nu_{n j}$. If $j$ is odd, this means $\nu_{n j} = 0$.}})
\beq \PIsot(u_0) \neq \emptyset \quad \Rightarrow \quad \nu_{n j} = 0~\mbox{for all $j$ odd}. \label{rec3} \feq
Independently of any assumption on $\PIsot(u_0)$, let us define
\beq \theta_n := \mbox{convergence radius of the series $\sum_{j=0}^{\infty} \nu_{n j} t^j$}
= \liminf_{j \vain + \infty} |\nu_{n j}|^{-1/j}~. \label{rooteta} \feq
Let us relate these objects to the convergence radius $\tau_n$ in \rref{taun}, to
the solution $u \in C^\om((-\Ti,T), \Aa)$
and to its squared $H^n$ norm.
We claim that
\beq \tau_n \leqs \theta_n~\mbox{and}~\| u(t) \|^2_n = \sum_{j=0}^{+\infty}
\nu_{n j} t^j~\mbox{for $t \in (-\Ti \wedge \theta_n, T \wedge \theta_n)$}
\label{claim0} \feq
(with $-\Ti \wedge \theta_n$ the opposite of $\Ti \wedge \theta_n$).
In fact:
the expansion $u(t) =  \sum_{j=0}^{+\infty} u_j t^j$, converging
in $\Hm{n}$ for $t \in (-\tau_n, \tau_n)$, implies
$\tau_n \leqs \theta_n$ and $\| u(t) \|^2_n = \sum_{j=0}^{+\infty}
\nu_{n j} t^j$ for $t \in (-\tau_n, \tau_n)$. Moreover the functions $t \in (-\Ti,T) \mapsto \| u(t) \|^2_n$
and $t \in (-\theta_n,\theta_n) \mapsto \sum_{j=0}^{+\infty}
\nu_{n j} t^j$ are analytic and coincide on $(-\tau_n, \tau_n)$,
so they coincide everywhere on the intersections of their domains,
which is $(-\Ti \wedge \theta_n, T \wedge \theta_n)$.
We now add to \rref{claim0} a stronger claim:
\beq \mbox{if $n > \dd{5 \over 2}$},~~
\tau_n \leqs \theta_n \leqs \Ti \wedge T~\mbox{and}~\| u(t) \|^2_n = \sum_{j=0}^{+\infty}
\nu_{n j} t^j~\mbox{for $t \in (-\theta_n, \theta_n)$}. \label{claim} \feq
Let us prove this claim, assuming for
example that $\Ti \wedge T = T$. If it were $T < \theta_n$ we would infer
$\lim_{t \vain T^{-}} \| u(t) \|^2_n = \lim_{t \vain T^{-}} \sum_{j=0}^{+\infty}
\nu_{n j} t^j = \sum_{j=0}^{+\infty}
\nu_{n j} T^j < + \infty$ (the first equality would hold
due to \rref{claim0} and $T \wedge \theta_n = T$; the subsequent two relations would hold because
$T$ would be inside the convergence interval of the series). On the other hand,
since $n > 5/2$, the conclusion that $\lim_{t \vain T^{-}} \| u(t) \|^2_n$ exists finite
would contradict \rref{blowhn}.
\section{Power series for the Euler equation in a paper of Behr, Ne$\boma{\check{\mbox{c}}}$as and  Wu}
\label{poneca}
In the paper \cite{Nec} mentioned above, the authors considered the power
series \rref{espa} for the Euler equation on $\Tt$, with an initial datum
$u_0 \in \Aa$ given by
\parn
\vbox{
\beq u_0(x) = \sum_{k = \pm a, \pm b, \pm c} u_{0 k} e^{i k \sc x}~, \label{unec} \feq
$$ a := (1,1,0),~~ b := (1,0,1),~~ c := (0,1,1)~; $$
$$ u_{0, \pm a} := (1,-1,0)~, \quad u_{0, \pm b}
:= (1,0,-1)~, \quad u_{0, \pm c} := (0,1,-1)~. $$
}
Like $u_0$, all the subsequent terms $u_j$ are Fourier
polynomials with rational coefficients ({\footnote{For a more precise statement
on these cofficients see our discussion of the datum $u_0$ in the next section
and, in particular, Eq. \rref{eqprec}.}).
Using rules equivalent to \rref{recur} \rref{repp}, the terms
$u_j$ were determined in \cite{Nec} by computer algebra,
for $j =1,2,...,35$. Computations were done with Mathematica for
$j=1,...,10$, and with a $C\scriptscriptstyle{++}$ program for $j=11,...,35$
(in the later case, approximating the rational coefficients
with finite precision decimal numbers). After determining
the $u_j$'s, the authors fixed their attention on the partial sums
$$ u^{(N)}(t) := \sum_{j=0}^N u_j t^j~, $$
whose $N \vain + \infty$ limit gives the solution
$u(t)$ of the Euler Cauchy problem, for all $t$ such that the series converges.
The previously mentioned computation of the $u_j$'s made available these partial
sums for $N=0,1,....,35$; the authors of \cite{Nec} computed the (squared) Sobolev norm
$$ \| u^{(N)}(t) \|^2_3  = \sum_{k \in \Zt} |k|^{6} |u^{(N)}_k(t)|^2 $$
for the above values of $N$, and several values of $t$.
Their main results were the following:
\parn
(i) Setting $t=0.32$, and analyzing the behavior of $\| u^{(N)}(0.32) \|_3$ for $N$
from $0$ to $35$, the authors found evidence
that $\| u^{(N)}(0.32) \|_3$ should approach a finite limit for $N \vain + \infty$. \parn
(ii) Setting $t=0.35$, the authors observed a rapid growth of
$\| u^{(N)}(0.35) \|_3$ for $N$ ranging from $0$ to $35$, a fact
suggesting that $\lim_{N \vain + \infty} \| u^{(N)}(0.35) \|_3 = +\infty$. \parn
(iii) A behavior as in (ii) was found to occur for slightly higher values of
$t$ (even though the authors suspected some rounding error to appear
for $t > 0.35$).
\par
The above results suggest that the series $\sum_{j=0}^{+\infty} u_j t^j$
has a finite convergence radius $\tau_3$ in $H^3$, with $\tau_3 \in (0.32, 0.35)$.
\par
Let us discuss this outcome from the viewpoint of the present paper,
denoting with $u$ the maximal $\Aa$-solution of the Cauchy problem
with this datum and recalling that this coincides with the maximal $H^3$-solution.
The datum $u_0$ possesses pseudo-symmetries (to be described in the next section);
therefore, $u$ has a time symmetric domain
$(-T,T)$ (in \cite{Nec} this fact was not explicitly declared, but probably
regarded as self-evident). According to our Eq. \rref{claimm}, it is
\beq \tau_3 \leqs T~; \feq
in principle, it could be $T = +\infty$.
In spite of this, the authors of \cite{Nec} spoke of a blow-up
at $\tau_3$.
\par
In the next two sections we present our computations
on the power series for the Behr-Ne$\boma{\check{\mbox{c}}}$as-Wu initial
datum, with our interpretation of the results.
Even though these calculations confirm the
''experimental'' outcomes (i)-(iii) of \cite{Nec},
we give evidence that the solution $u$ of the Euler
equation does not blow up close to $\tau_3$; on the contrary,
computing the power series for $\| u(t) \|^2_3$ up
the available order we obtain strong evidence
that such a power series has a convergence radius
$\theta_3$ such that $0.47 < \theta_3 < 0.50$,
which implies for the time $T$ of existence
of $u$ the bound $T \geqs \theta_3 > 0.47$.
By a subsequent analysis relying on the technique of the \Pade approximants,
we show that a blow-up of $u(t)$
might happen at a time larger than $0.48$: more precisely,
these computations give a somehow weak indication
that $T$ might be finite, with $0.56 < T < 0.73$.
\section{Our approach to the power series of Behr, Ne$\boma{\check{\mbox{c}}}$as and  Wu}
\label{ourap}
Let us denote again with $u_0$ the datum \rref{unec} and consider its iterates
$u_{j}$ ($j=1,2,...$), with the corresponding power series;
like $u_0$, all the iterates $u_j$ are Fourier polynomials
with rational coefficients. Throughout the section, $u$ is the
maximal $\Aa$-solution of the Euler equation with datum $u_0$.
\salto
\textbf{A closer analysis of the Behr-Ne$\check{\mbox{c}}$as-Wu initial datum:
symmetry properties.}
The symmetry group $\Isot(u_0)$ and the pseudo-symmetry space
$\PIsot(u_0)$ (Eqs. \rref{isot} \rref{pisot}) can be explicitly computed.
For the first one, we find
\parn
\vbox{
\beq \Isot(u_0) = \{ (\uno, 0), (\uno, \imath_2),
(A, a_1), (A, a_2), (B, a_1), (B, a_2), \label{isotnec} \feq
$$ (C, c_1), (C, c_2), (D, c_1), (D, c_2), (E, 0), (E, \imath_2) \} $$
}
\parn
where $\uno$ is the $3 \times 3$ identity matrix, and
\beq A := \left( \barray{ccc} 0 & 1 & 0 \\ 0 & 0 & 1 \\ 1 & 0 & 0 \farray \right),
\quad B := \left( \barray{ccc} 0 & -1 & 0 \\ -1 & 0 & 0 \\ 0 & 0 & -1 \farray \right)  \label{ae} \feq
$$ C := \left( \barray{ccc} 0 & 0 & 1 \\ 1 & 0 & 0 \\ 0 & 1 & 0 \farray \right),
\quad D := \left( \barray{ccc} -1 & 0 & 0 \\ 0 & 0 & -1 \\ 0 & -1 & 0 \farray \right)
\quad E := \left( \barray{ccc} 0 & 0 & -1 \\ 0 & -1 & 0 \\ -1 & 0 & 0 \farray \right)~; $$
furthermore, $\imath_2, a_1$, etc., are the following elements of $\Tt$:
\beq
\imath_2 := (\pi, \pi, \pi)~,~~
a_1 := (0,0,\pi)~,~~ a_2 := (\pi, \pi,0)~, \label{a1e2} \feq
$$ c_1 := (\pi,0,0)~,~~ c_2 := (0,\pi, \pi)~ $$
(of course, in the above $\pi$ is short for $\pi$ mod. $2 \pi \interi$).
Let us fix the attention on the reduced symmetry subgroup $\Isott(u_0) = \{\uno, A,B,C,D,E \}$; it
is readily checked that
\beq A^3 = \uno, B^2 = \uno, (B A)^2 = \uno \label{eqgen} \feq
$$ C = A^2, \quad D = A B~, \quad E = A^2 B~. $$
So, $\Isott(u_0)$ has two generators $A, B$; the first line
in \rref{eqgen} gives a presentation of this group in terms
of generators and relations, while the second line expresses
the other elements in terms of $A, B$. Using Eq. \rref{eqgen},
one recognizes a group isomorphism
\beq \Isott(u_0) \simeq \textbf{D}_3 \feq
where the right-hand side indicates the dihedral group of order 3,
formed by the symmetries of an equilateral triangle
({\footnote{For any integer $n \in \{3,4,....\}$, one denotes
with $\textbf{D}_n$ the dihedral group of order $n$;
this is formed by the orthogonal transformations of
the Euclidean plane $\reali^2$ into itself which preserve
a regular polygon with $n$ sides, centered at the origin.
Denoting with $id$ the identity map,
with $a$ the rotation of an angle $2 \pi/n$
and with $b$ the reflection about anyone of the $n$ symmetry axes of the
polygon, one finds that $a, b$ are generators of $\textbf{D}_n$
and fulfill the relations $a^n = id$, $b^2 = id$, $(b a)^2 = id$.
The elements of $\textbf{D}_n$ are $2 n$, and coincide
with $id, a,a^2,...,a^{n-1}$, $b$, $a b$, $a^2 b$, ..., $a^{n-1} b$.}}).
\parn
\vbox{
\vskip 0.08cm
\noindent
Now we consider the full group $\Isot(u_0)$
(with the product \rref{prcirc}). It is easy to check that
\beq
(A, a_1)^6 = (\uno,0)~,~~(B, a_1)^2 = (\uno,0)~, ((B, a_1)(A, a_1))^2 = (\uno,0)~,
\label{eqqgen} \feq
$$ (A, a_1)^2 = (C, c_2)~,~~(A, a_1)^3 = (\uno, \imath_2)~,~~
(A, a_1)^4 = (A, a_2)~,~~(A, a_1)^5 = (C, c_1)~, $$
$$ (A, a_1) (B, a_1) = (D, c_2)~,~~(A, a_1)^2 (B, a_1) = (E, \imath_2)~,~~
(A, a_1)^3 (B,a_1) = (B, a_2)~, $$
$$ (A, a_1)^4 (B, a_1) = (D, c_1)~,~~
(A, a_1)^5 (B, a_1) = (E, 0)~. $$
}
\parn
So, $\Isot(u_0)$ has two generators $(A, a_1)$ and  $(B, a_1)$; the first line
in \rref{eqqgen} gives a presentation of this group in terms
of generators and relations, and the subsequent lines express
the other elements in terms of the generators. One recognizes a group isomorphism
\beq \Isot(u_0) \simeq \textbf{D}_6 \feq
where the right-hand side indicates the dihedral group of order 6,
formed by the symmetries of a hexagon (see the previous footnote).
\par
Let us pass to the pseudo-symmetry space $\PIsot(u_0)$. One
readily checks that this contains $(-\uno,0)$
(inducing the space reflection $\E(-\uno,0) : x \in \Tt \mapsto -x$).
From here and from the general result
\rref{pseudorifl}, one obtains
\parn
\vbox{
\beq \PIsot(u_0) =  \Isot(u_0) (-\uno,0)  = \{ (-\uno, 0), (-\uno, \imath_2),
(-A, a_1), (-A, a_2),  \label{pisotnec} \feq
$$ (-B, a_1), (-B, a_2), (-C, c_1), (-C, c_2), (-D, c_1), (-D, c_2), (-E, 0), (-E, \imath_2) \}~, $$
}
with $A,B,...$ and $\imath_2,a_1,...$ as in Eqs. \rref{ae} \rref{a1e2}.
\salto
\textbf{Some consequences of the previous symmetry results.}
(i) What we have stated in Section \ref{power} for an arbitrary initial datum
holds, in particular, for the present datum $u_0$:
the symmetries or pseudo-symmetries of $u_0$
can be used to speed up the computation of the Fourier
components of any iterate $u_j$. More precisely, if we know
the Fourier component $u_{j, k}$ for some $k$,
using Eqs. \rref{equal2} \rref{pequal2} we readily
obtain the components $u_{j, S k}$ for all $S \in \Isott(u_0)
\cup \PIsott(u_0)$. \parn
(ii) As already noted, the pseudo-symmetry space $\PIsot(u_0)$ contains
$(-\uno, 0)$, corresponding to the space reflection.
In terms of Fourier coefficients,
the relation \rref{pequal2} with $(S, a) = (-\uno, 0)$ takes
the form $u_{j, - k} = (-1)^{j} u_{j, k}$ for $j=0,1,2,...$ and $k \in \Zt$.
On the other hand, any iterate $u_{j}$ is a real vector field, thus
$u_{j, - k} = \overline{u_{j, k}}$; in conclusion
$\overline{u_{j, k}} = (-1)^{j} u_{j, k}$, which indicates that
$u_{j, k}$ is real for $j$ even, and imaginary for $j$ odd. Taking into
account that the coefficients $u_{j, k}$ are rational in any case, we conclude
the following for each $k \in \Zt$:
\beq u_{j, k} \in \razionali^3~\mbox{for $j=0,2,4,...$}~; \qquad
u_{j, k} \in i \razionali^3~\mbox{for $j=1,3,5,...$}~. \label{eqprec} \feq
(iii) In the sequel we are often interested in
the partial sums
$u^{(N)}(t) := \sum_{j=0}^N u_j t^j$
and in their norms $\| u^{(N)}(t) \|_n$, especially for $n=3$.
Since $\PIsot(u_0) \neq \emptyset$, as in
\rref{norequal} we have $\| u_N(t) \|_n = \| u_N(-t) \|_n$. \parn
(iv) Independently of any convergence consideration about
the power series $\sum_{j=0}^{+\infty} u_j t^j$, the
result $\PIsot(u_0) \neq \emptyset$ also ensures
that the (maximal $\Aa$-) solution $u$ of the Euler equation
with datum $u_0$ has a symmetric domain $(-T,T)$
(recall Eq. \rref{rs22}).
\salto
\textbf{Describing our computations.}
We have considered again the power series \rref{espa}
for the datum $u_0$; to deal with this series
we have written a program in
Python, using the package
gmpy \cite{gmpy} for fast arithmetics on rational numbers.
This program implements Eq. \rref{repp} for $\PPP$ and the recursion rule
\rref{recur}; moreover, it takes into account
the dihedral symmetries (and pseudosymmetries) of $u_0$ to speed up computations.
The program has been run to compute
the terms $u_j$ for $j=1,...,52$
({\footnote{To test the reliability of this program,
the calculation of some of the $u_j's$ has been
checked in two independent ways. These checks have been
done by means of other two programs, which
implement Eqs. \rref{repp} \rref{recur} accepting
as an initial datum $u_0$ any Fourier polynomial;
these do not refer to any symmetry property of $u_0$. The first
of these programs, written in Mathematica, has
been used to compute the $u_j$'s up to order $j=13$;
the second program, written in Python, has been used for a calculation
up to $j=43$.}}).
Calculations have been
performed on a PC with an Intel Core i7 CPU 860 at 2.8GHz
and an 8GB RAM.
The CPU time for $u_j$ has been, for example:
$1$ second for $j=10$, one minute for $j=20$,
half an hour for $j=30$,
$7$ hours for $j=40$ and $85$ hours for $j=52$.
Differently from \cite{Nec}, for \textsl{all} orders up to $j=52$ the Fourier coefficients
$u_{j, k}$ of $u_j$ have been represented as
elements of $\razionali^3$ or
$i \razionali^3$; so, no rounding errors
related to finite precision arithmetics have been
introduced in the calculation of the power series. \par
From the $u_j$'s one determines the squared norms
$\| u_j \|^2_3 = \sum_{k \in \Zt} |k|^{6} |u_{j, k}|^2$,
the partial sums
$u^{(N)}(t) := \sum_{j=0}^N u_j t^j$
and their squared norms $\| u^{(N)}(t) \|^2_3$  $= \sum_{k \in \Zt} $ $|k|^{6} |u^{(N)}_k(t)|^2$
($N=1,...,52)$.
Each $\| u_{j} \|^2_3$ is a rational number
and $\| u^{(N)}(t) \|^2_3$ is a
polynomial of order $2 N$ in $t$, with rational coefficients,
containing only even powers of $t$; furthermore, the coefficients
of $t^0$ and $t^{2 N}$ in $\| u^{(N)}(t) \|^2_3$ are $\| u_0 \|^2_3$ and
$\| u_{N} \|^2_3$, respectively. \par
Our computations of the above norms, up to $j=52$ or $N=52$, have been
done using the previously mentioned Python program. These calculations
have been relatively quick: for example, the computation of
$\| u^{(52)}(t) \|^2_3$ has required a CPU time of about $3$ hours.
As first examples of our results, we report the following ones:
\beq \| u_0 \|^2_3 = 96, \quad \| u_1 \|^2_3 = 6912, \quad \| u_2 \|^2_3 = 45440, \feq
$$ \| u_3 \|^2_3 = {3695360 \over 9}, \quad \| u_4 \|^2_3 = {1366793248 \over 675},
\quad \| u_5 \|^2_3 = {2243123779689032 \over 186046875}~. $$
$\| u_{52} \|^2_3$ is a ratio of integers where the numerator and
the denominator have $19515$ and $19463$ digits, respectively.
Table 1 reports $\| u_{j} \|^2_{3}$ for
$j=0,...,52$, in the $16$ digits decimal representation.
\vskip 0.4cm
\centerline{\textbf{Table 1}. The squared norms $\| u_{j} \|^2_3$.}
\vskip 0.2cm
\begin{table}[h]
\centering
\subtable{
\begin{tabular}{|l|l|}
\hline
$j$ & $\| u_j \|^2_3$ \\
\hline
$1$ &  $6912$ \\
$2$ &  $45440$ \\
$3$ &  $4.105955555555556 \times 10^5$ \\
$4$ &  $2.024878885925926 \times 10^{6}$ \\
$5$ &  $1.205676676745595 \times 10^{7}$ \\
$6$ &  $8.452219877103332 \times 10^{7}$ \\
$7$ &  $6.152775603322622 \times 10^{8}$ \\
$8$ &  $4.791192836997696 \times 10^{9}$ \\
$9$ &  $3.628869598772102 \times 10^{10}$ \\
$10$ &  $2.825486371143428 \times 10^{11}$ \\
$11$ &  $2.228507964437443 \times 10^{12}$ \\
$12$ &  $1.821213808657725 \times 10^{13}$ \\
$13$ &  $1.539790191793044 \times 10^{14}$ \\
$14$ &  $1.341372343677860 \times 10^{15}$ \\
$15$ &  $1.190159209731028 \times 10^{16}$ \\
$16$ &  $1.066432595016119 \times 10^{17}$ \\
$17$ &  $9.598519025230687 \times 10^{17}$ \\
$18$ &  $8.662788463495777 \times 10^{18}$ \\
$19$ &  $7.840631870939454 \times 10^{19}$ \\
$20$ &  $7.122921654632158 \times 10^{20}$ \\
$21$ &  $6.499436510134908 \times 10^{21}$ \\
$22$ &  $5.957837347113741\times 10^{22}$ \\
$23$ &  $5.485035371335649\times 10^{23}$ \\
$24$ &  $5.068929708200902\times 10^{24}$ \\
$25$ &  $4.699401376031744\times 10^{25}$ \\
$26$ &  $4.368534165204974\times 10^{26}$ \\
\hline
\end{tabular}
}
\hspace{-0.6cm}
\subtable{
\begin{tabular}{|l|l|}
\hline
$j$ & $\| u_j \|^2_3$ \\
\hline
$27$ & $4.070323867244879 \times 10^{27}$ \\
$28$ & $3.800202819232687 \times 10^{28}$ \\
$29$ & $3.554589555246873 \times 10^{29}$ \\
$30$ & $3.330557264153261 \times 10^{30}$ \\
$31$ & $3.125627141295364 \times 10^{31}$ \\
$32$ & $2.937654907691943 \times 10^{32}$ \\
$33$ & $2.764771414352126 \times 10^{33}$ \\
$34$ & $2.605347861791808 \times 10^{34}$ \\
$35$ & $2.457968790658826 \times 10^{35}$ \\
$36$ & $2.321406184470901 \times 10^{36}$ \\
$37$ & $2.194593722846032 \times 10^{37}$  \\
$38$ & $2.076602420620089 \times 10^{38}$ \\
$39$ & $1.966618988613002 \times 10^{39}$  \\
$40$ & $1.863927582086700 \times 10^{40}$ \\
$41$ & $1.767894900465337 \times 10^{41} $  \\
$42$ & $1.677958174980847 \times 10^{42}$ \\
$43$ & $1.593615440091581 \times 10^{43}$ \\
$44$ & $1.514417532673484 \times 10^{44}$ \\
$45$ & $1.439961389630372 \times 10^{45}$ \\
$46$ & $1.369884345517744 \times 10^{46}$ \\
$47$ & $1.303859232337703 \times 10^{47}$ \\
$48$ & $1.241590147950303 \times 10^{48}$ \\
$49$ & $1.182808795435820 \times 10^{49}$ \\
$50$ & $1.127271314453561 \times 10^{50}$ \\
$51$ & $1.074755536205362 \times 10^{51}$ \\
$52$ & $1.025058601409640 \times 10^{52}$ \\
\hline
\end{tabular}
}
\end{table}
\vfill \eject \noindent
{~}
\parn
\vbox{
\noindent
Let us pass to the squared norms $\| u^{(N)}(t) \|^2_3$. As an example, the result for
$N=5$ is
\beq \| u^{(5)}(t) \|^2_3 = 96 + 6656 \, t^2 + {258304 \over 9} \, t^4 + {104566912 \over 525} \, t^6 \feq
$$ - {9513575648 \over 70875} \, t^8 + {2243123779689032 \over 186046875} \, t^{10}~.$$
}
There is no room to report here the results obtained for all the other values of $N$, especially
in the rational form for the coefficients. However, we can write some of them in the $16$ digits
precision; in particular,
\parn
\vbox{
\beq \| u^{(52)}(t) \|^2_3 \label{sedici} \feq
$$ = 96 + 6656 \, t^2 + 2.870044444444444\times 10^4 \, t^4 + 1.993359937918871\times 10^5 \, t^6 $$
$$ + 1.058054454761424\times 10^5 \, t^8 + 1.781444415306641 \times 10^6 \, t^{10}
+  2.740017914111055 \times 10^6 \, t^{12} $$
$$ - 7.321985472578865 \times 10^6 \, t^{14}
+ 4.183410651491110 \times 10^{6} \, t^{16} + 1.457483700816015 \times 10^8 \, t^{18} $$
 $$ - 1.768517246168822 \times 10^8 \, t^{20} + 4.196205149715839 \times 10^8 \, t^{22}
 + 3.648789154816725 \times 10^9 \, t^{24} $$
 $$  - 2.178830191383206 \times 10^{10} \, t^{26}
 - 1.394064522752687 \times 10^{10} \, t^{28} + 2.954202883502504 \times 10^{11} \, t^{30} $$
 $$ + 1.283692616423054 \times 10^{11} \, t^{32} - 4.543575106022102 \times 10^{12} \, t^{34}
 +  4.789569007452901 \times 10^{12} \, t^{36} $$
 $$ + 2.830635227431622 \times 10^{13} \, t^{38}
 +  4.470168139346678 \times 10^{13} \, t^{40} - 6.910532995061547 \times 10^{14}\,t^{42} $$
 $$ + 1.457019276470951 \times 10^{14}\,t^{44} + 9.053007124662626 \times 10^{15}\,t^{46}
 - 8.939780851014422  \times 10^{15} \,t^{48} $$
 $$ - 1.019952729404346 \times 10^{17}\,t^{50} + 1.137772938577812 \times 10^{17}\,t^{52}
 + 1.644161010427522 \times 10^{18} \,t^{54} $$
 $$ - 4.571936581656874 \times 10^{18}\,t^{56} - 3.140936865806385 \times 10^{19} \,t^{58}
 + 2.408085513008218 \times 10^{21}\,t^{60} $$
 $$ - 1.107900217253947 \times 10^{23}\,t^{62}
 + 4.186064092726056 \times 10^{24}\,t^{64} - 1.674853723772203 \times 10^{26}\,t^{66} $$
 $$ + 6.911508987260593 \times 10^{27}\,t^{68} - 2.698282390313396 \times 10^{29}\,t^{70}
 + 9.951375797771149 \times 10^{30}\,t^{72} $$
$$ - 3.558771163372845 \times 10^{32}\,t^{74}
+ 1.232107326257251 \times 10^{34}\,t^{76} - 4.045044388392564 \times 10^{35}\,t^{78} $$
$$ + 1.242344004413423 \times 10^{37}\,t^{80} - 3.561397641466941 \times 10^{38}\,t^{82}
 + 9.520206481050174 \times 10^{39} \,t^{84} $$
 $$ - 2.357932432543021 \times 10^{41}\,t^{86}
 + 5.354229494719748 \times 10^{42}\,t^{88}
 - 1.103667607665446 \times 10^{44} \,t^{90} $$
 $$ + 2.052382635232918 \times 10^{45}\,t^{92}
 - 3.436006560519912 \times 10^{46}\,t^{94}
 + 5.184487278969682 \times 10^{47} \,t^{96} $$
 $$ - 7.072466985323957 \times 10^{48}\,t^{98}
 + 8.759614973466463 \times 10^{49}\,t^{100}
 - 9.896987665647683 \times 10^{50}\,t^{102} $$
$$ + 1.025058601409640 \times 10^{52}\,t^{104} ~. $$
}
The rest of the paper reports a number of facts stemming from our computations, with the interpretation that we
suggest for them.
\vfill \eject \noindent
{~}
\vskip 1cm \noindent
\textbf{Verification of the outcomes of \cite{Nec} on
$\boma{\| u^{(N)}(t) \|^2_3}$.}
Our computations based on the systematic use of rational
numbers have given essentially the same results as in \cite{Nec}
about $\| u^{(N)}(t) \|^2_3$ as a function of $N$,
in the two cases $t=0.32$ and $t=0.35$. So,
$\| u^{(N)}(0.32) \|^2_3$ seems to approach a
limit value for large $N$, while $\| u^{(N)}(0.35) \|^2_3$
grows rapidly with $N$; our use of rational coefficients
ensures that such a rapid growth
is not due to cumulative rounding errors.
In Figures 1-2, we report $\| u^{(N)}(t) \|^2_3$
as a function of $N \in \{0,...,52\}$, in the two cases
$t=0.32$ and $t=0.35$; these figures are very similar to the
ones at the bottom of pages 235 and 236 of \cite{Nec}, respectively
(but comparison requires a rescaling, since the $H^3$ norm
employed in \cite{Nec} differs from ours by a constant factor).
\par
We agree with \cite{Nec} in interpreting these results
as indications that the power series for this initial datum has a finite
$H^3$-convergence radius $\tau_3$, with $\tau_3 \in (0.32, 0.35)$.
\vskip 1cm \noindent
\begin{figure}[h]
\parbox{3in}{
\includegraphics[
height=2.0in,
width=2.8in
]%
{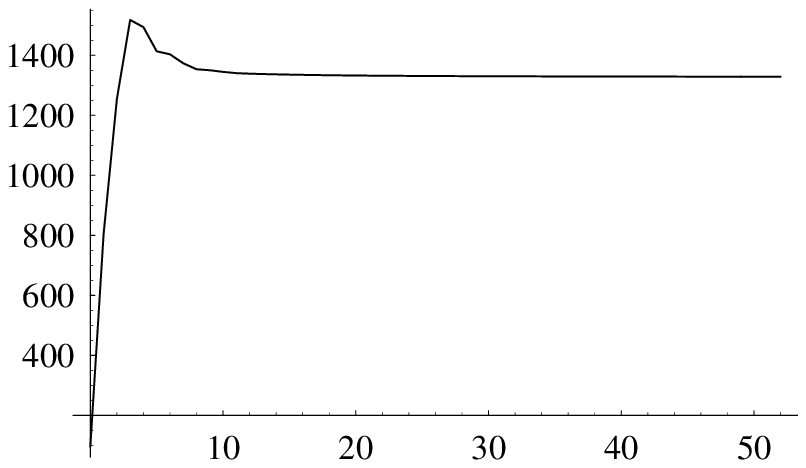}%
\par
{\footnotesize{
{\textbf{Figure 1.~} $\| u^{(N)}(0.32) \|^2_3$
as a function of \\ $N \in \{0,1,...,52\}$.
}}
\par}
\label{f1}
}
\hskip 0.4cm
\parbox{3in}{
\includegraphics[
height=2.0in,
width=2.8in
]%
{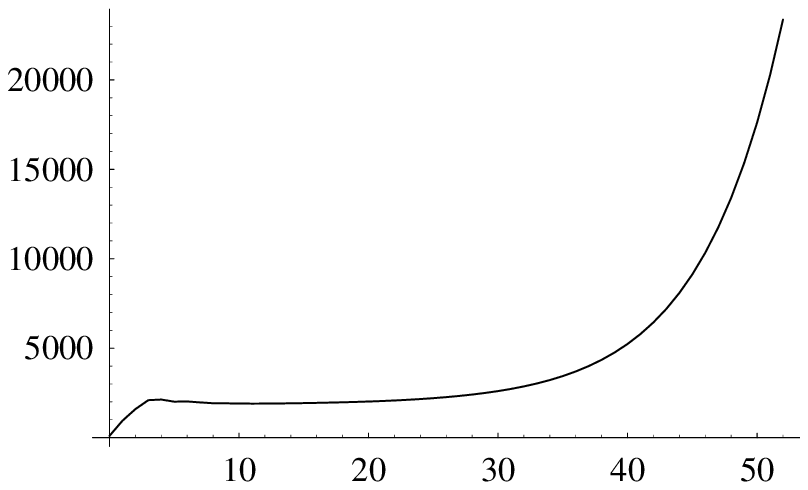}%
\par
{\footnotesize{
{\textbf{Figure 2.~} $\| u^{(N)}(0.35) \|^2_3$
as a function of \\ $N \in \{0,1,...,52\}$.
}}
\par}
\label{f2}
}
\end{figure}
\vfill \eject \noindent
\textbf{Further evidence on the $\boma{H^3}$-convergence radius of
the power series.} This comes from the root test \rref{roote} for $n=3$:
\beq \tau_{3} = \liminf_{j \vain + \infty} \| u_j \|_3^{-1/j}~. \feq
Figure 3 represents $\| u_j \|_3^{-1/j}$
as a function of $j \in \{1,...,52\}$.
For $j=36,38,...,52$ we have a very good interpolation
\beq \| u_j \|_3^{-1/j} \simeq 0.32158 - \left(1.20125 \over j \right)^{1.38458}~, \label{vegood} \feq
(obtained assuming for the interpolant the form $A - (B/j)^c$, and
applying the least squares criterion).
The right-hand side of \rref{vegood} approximates $\| u_j \|_3^{-1/j}$
with a mean quadratic error $< 10^{-5}$ (averaging,
as indicated, for $j=36,38,...,52$; if we
average over the larger range
$j=16,18,...,52$, the mean quadratic error is $< 10^{-4}$). \par
Assuming that \rref{vegood} approximates
$\| u_j \|_3^{-1/j}$ with a similar precision for arbitrarily large $j$, but keeping prudentially
only two digits in our final estimate, we conclude with an estimate
\beq 0.32 < \tau_3 < 0.33~. \label{estau3} \feq
\begin{figure}[h]
\hspace{3cm}
\parbox{3in}{
\includegraphics[
height=2.0in,
width=2.8in
]%
{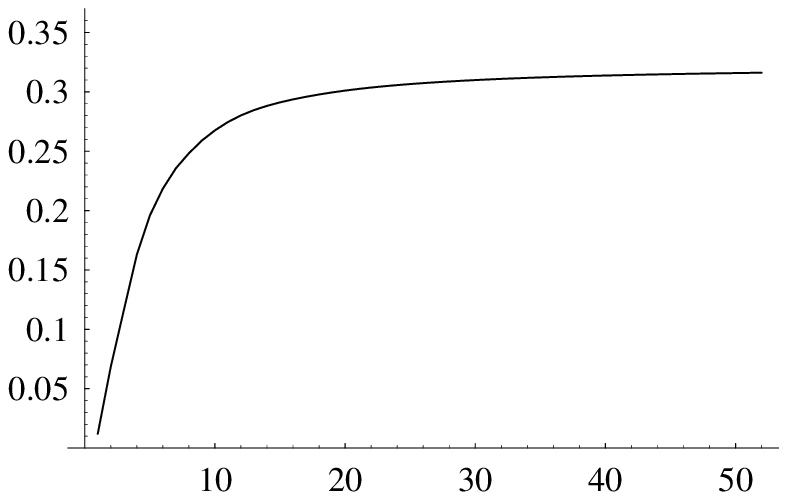}%
\par
{\footnotesize{
{\textbf{Figure 3.} \hbox{$\| u_j \|^{-1/j}_3$
as a function of $N \in \{1,...,52\}$.}
}}
\par}
\label{f3}
}
\end{figure}
\vfill \eject \noindent
\textbf{Reminder estimates for the series expansion of $\boma{u(t)}$ in $\boma{\Hm{3}}$.}
Let $N \in \{0,1,2,...\}$; of course
\beq u(t) - u^{(N)}(t) = \sum_{j=N+1}^{+\infty} u_j t^j \qquad \mbox{for $t \in (-\tau_3,\tau_3)$}~; \feq
this implies
\beq \| u(t) - u^{(N)}(t) \|_3 \leqs \sum_{j=N+1}^{+\infty} \| u_j \|_3 |t|^j \qquad
\mbox{for $t \in (-\tau_3,\tau_3)$}~. \label{rem2} \feq
To go on, we need a guess on the behavior of the norms $\| u_j \|_3$.
To this purpose, let us consider the sequence
\beq  \mu_{3 j}:= 0.32^j \, \| u_j \|_3 \qquad  (j=0,1,2...)~, \feq
recalling that $0.32$ is the lower bound
for $\tau_3$ in \rref{estau3}. From the norms available up to $j=52$, we can check that $(\mu_{3 j})$ is decreasing
while $j$ ranges in $\{1,3,...,52\}$; by extrapolation, let us assume that
$(\mu_{3 j})$ is decreasing on the infinite set $\{1,2,....\}$. So, $\mu_{3 j} \leqs \mu_{3 N}$
for integer $j \geqs N \geqs 1$, i.e.,
\beq \| u_j \|_3 \leqs {\mu_{3 N} \over 0.32^j} \qquad \mbox{for $j \geqs N \geqs 1$}~. \feq
For $t \in (-0.32, 0.32)$, inserting this inequality into \rref{rem2}
we get
$\| u(t) - u^{(N)}(t) \|_3 \leqs \mu_{3 N} \sum_{j=N+1}^{+\infty} |t/0.32|^j =
\mu_{3 N} |t/0.32|^{N+1} \sum_{j=0}^{+\infty} |t/0.32|^j$, i.e.,
\beq \| u(t) - u^{(N)}(t) \|_3 \leqs \mu_{3 N} {|t/0.32|^{N+1} \over 1 - |t/0.32|}
~~\mbox{for $t \in (-0.32, 0.32)$, $N \in \{1,2,3,...\}$}. \label{rem3} \feq
Of course, this is a conjecture based on the previous extrapolation. For the practical
application of the reminder estimate \rref{rem3}, we mention that (rounding up
from above)
\beq \mu_{3 \, 5} =  11.7,~ \mu_{3 \, 10} = 5.99,~ \mu_{3 \, 20} = 3.39,~ \feq
$$ \mu_{3 \, 30} = 2.61,~ \mu_{3 \, 40} = 2.20, ~\mu_{3 \, 52} = 1.88 \,. $$
\vskip 0.2cm\noindent
\textbf{No blow-up at $\boma{\tau_3}$.} After accumulating indications
that the Taylor series for $u(t)$ has an $H^3$-convergence radius
$\tau_3 \in (0.32,0.33)$, in the rest of the section we will present
evidence that $u(t)$ does \textsl{not} blow up at $t=\tau_3$.
\vskip 0.2cm \noindent
\textbf{The power series for $\boma{\| u(t) \|^2_3}$;
an indication that $\boma{u(t)}$ exists
up to $\boma{t=0.47}$ at least.}
The results \rref{formal} \rref{rec3} with $n=3$ give
a formal series expansion
\beq \| \sum_{j=0}^{+\infty} u_j t^j \|^2_3 = \sum_{j=0}^{+\infty} \nu_{3 j} t^j,
\quad \nu_{3 j} := \sum_{\ell=0}^j \la u_\ell | u_{j - \ell} \ra_3 \in \reali,
\quad \nu_{3 j} = 0 \quad \mbox{for $j$ odd}~; \feq
the series $\sum_{j=0}^{+\infty} \nu_{3 j} t^j$ has a convergence radius
\beq \theta_3 = \liminf_{j \vain + \infty} |\nu_{3 j}|^{-1/j}~. \label{rooteta3} \feq
Recalling that $(-T,T)$ is the domain of the solution $u$, we know (from
\rref{claim}) that
\beq \tau_3 \leqs \theta_3 \leqs T~, \qquad \| u(t) \|^2_3 = \sum_{j=0}^{+\infty}
\nu_{3 j} t^j \qquad \mbox{for $t \in (-\theta_3, \theta_3)$}. \label{claim3} \feq
In the sequel, for $N=0,1,2,....$ we also consider the partial sums
\beq \nu^{(N)}_3(t) := \sum_{j=0}^{N} \nu_{3 j} t^j~. \feq
Of course,
$u^{(N)}(t) = \sum_{j=0}^N u_j t^j$ is such that
$u(t) = u^{(N)}(t) + O(t^{N+1})$ for $t \vain 0$; this implies
$\| u(t) \|^2_3 = \| u^{(N)}(t) \|^2_3 + O(t^{N+1})$, whence
\beq  \nu^{(N)}_3(t)  = {\| u^{(N)}(t) \|^2_3~\Big|_{\scriptscriptstyle{\mbox{$t^k \vain 0$ for $k > N$}}}}
\hspace{0.2cm}.
\label{der1} \feq
With this remark, the previuos computations of $\| u^{(N)}(t) \|^2_3$
up to $N=52$ also give the partial sums $\nu^{(N)}_3(t)$ for
$N=1,...,52$ or, equivalently, the coefficients $\nu_{3 j}$ for
$j=0,...,52$. For example,
\beq \nu_{3 0}  = 96, \quad \nu_{3 2} = 6656, \quad
\nu_{3 4} = {258304 \over 9}, \quad \nu_{3 6} = {2825587712 \over 14175}, \feq
$$ \nu_{3 8}  = {52545219363488 \over 496621125}, \quad
\nu_{3 \, 10}  = {10025320340466597351685768 \over 5627635784943046875}~; $$
$\nu_{3 \, 52}$ is a ratio of integers where the numerator and
the denominator have $2610$ and $2593$ digits, respectively.
\par
The $16$-digits representation of the coefficients
$\nu_{3 j}$ for all $j \in \{0,...,52\}$ can be obtained from
Eqs. \rref{sedici} \rref{der1}; more precisely,
\beq \nu_{3 j} = \mbox{coefficient of $t^j$ in \rref{sedici}, for $j=0,...,52$}~. \feq
\vfill \eject
From the above data, one can try to make predictions on the convergence
radius $\theta_3$ of the series $\sum_{j=0}^{+\infty} \nu_{3 j} t^j$.
In Figures 4-7 we report the partial sums
$\nu^{(N)}_3(t)$
as functions of $N \in \{0,...,52\}$, in the four cases
$t=0.45, 0.50, 0.55, 0.60$. For $t=0.45$, the function $N \mapsto
\nu^{(N)}_3(t)$ seems to approach a limit value
for large $N$. The situation is not clear for $t=0.50$,
due to the appearing of small oscillations; for
$t=0.55$ and $t=0.60$, the oscillations of $N \mapsto \nu^{(N)}_3(t)$
are large and their amplitude increases
with $N$. We regard these results as indicating
that $\sum_{j=0}^{+\infty} \nu_{3 j} t^j$ is
convergent for $t \leqs 0.45$ and not convergent
for $t \geqs 0.55$; in other words, for the convergence
radius we have a conjectural estimate
\beq 0.45 < \theta_3 <  0.55~. \label{pteta3} \feq
{~}
\begin{figure}[h]
\parbox{3in}{
\includegraphics[
height=2.0in,
width=2.8in
]%
{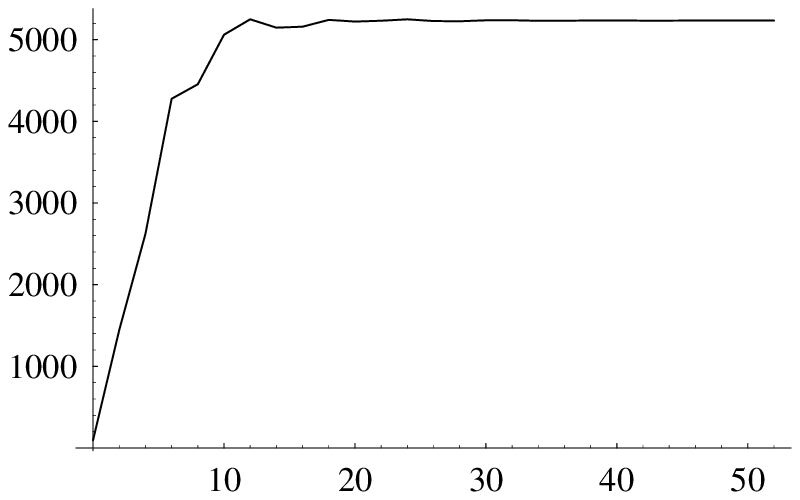}%
\par
{\footnotesize{
{\textbf{Figure 4.~} $\nu^{(N)}_3(0.45)$
as a function of \\ $N \in \{0,2,...,50,52\}$.
}}
\par}
\label{f45}
}
\hskip 0.4cm
\parbox{3in}{
\includegraphics[
height=2.0in,
width=2.8in
]%
{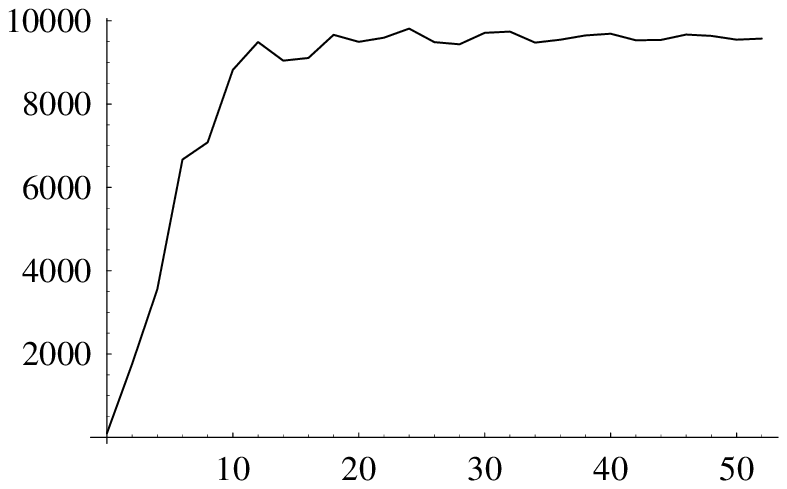}%
\par
{\footnotesize{
{\textbf{Figure 5.~} $\nu^{(N)}_3(0.50)$
as a function of \\ $N \in \{0,2,...,50,52\}$.
}}
\par}
\label{f50}
}
\end{figure}
\begin{figure}[h]
\parbox{3in}{
\includegraphics[
height=2.0in,
width=2.8in
]%
{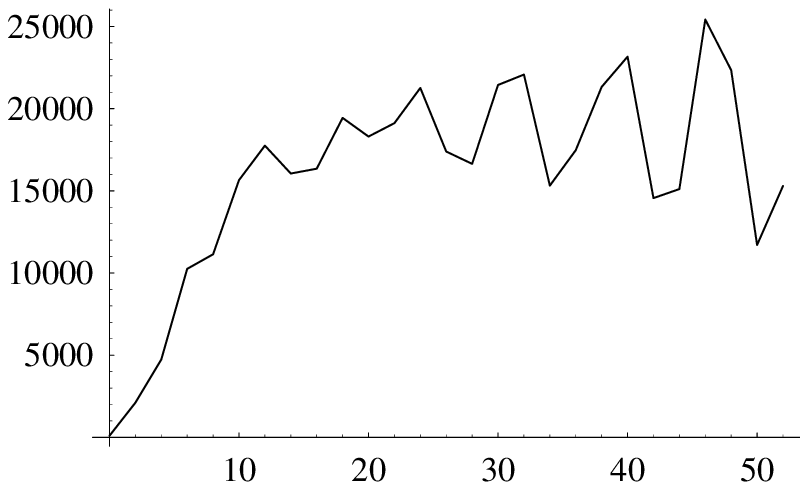}%
\par
{\footnotesize{
{\textbf{Figure 6.~} $\nu^{(N)}_3(0.55)$
as a function of \\ $N \in \{0,2,...,50,52\}$.
}}
\par}
\label{f55}
}
\hskip 0.4cm
\parbox{3in}{
\includegraphics[
height=2.0in,
width=2.8in
]%
{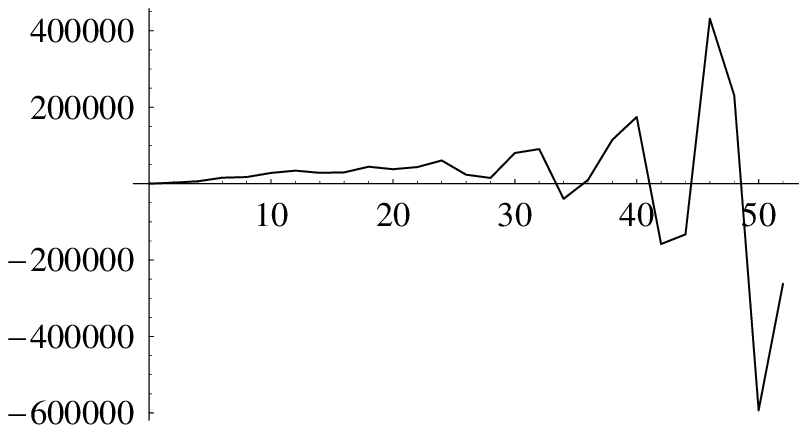}%
\par
{\footnotesize{
{\textbf{Figure 7.~} $\nu^{(N)}_3(0.60)$
as a function of \\ $N \in \{0,2,...,50,52\}$.
}}
\par}
\label{f60}
}
\end{figure}
\vfill \eject
Another way to estimate $\theta_3$ comes
from the root test \rref{rooteta3}. Figure 8 is a graph of
$|\nu_{3 j}|^{-1/j}$ as a function of $j \in \{2,4,...,50,52\}$.
For $j=36,38,...,52$, there is a fairly good interpolation
\beq |\nu_{3 j}|^{-1/j} \simeq 0.484 - \left({8.48 \over j}\right)^{2.19}
\label{nu3jin} \feq
(obtained assuming for the interpolant the form $A - (B/j)^c$, and
applying the least squares criterion); here,
the right-hand side approximates $|\nu_{3 j}|^{-1/j}$
with a mean quadratic error $< 0.01$
(let us repeat it, for $j$ between $36$ and $52$).
Assuming that the above interpolant behaves similarly for
all larger (even) $j$, and considering
$\theta_3 = \liminf_{j \vain + \infty} |\nu_{3 j}|^{-1/j}$ we are led
to use $0.484 \pm 0.01$ as upper and lower bounds for it; rounding up
to two digits we obtain the inequality
\beq 0.47 < \theta_3 < 0.50~, \label{teta3} \feq
which is compatible with \rref{pteta3}.
\begin{figure}[h]
\hspace{3cm}
\parbox{3in}{
\includegraphics[
height=2.0in,
width=2.8in
]%
{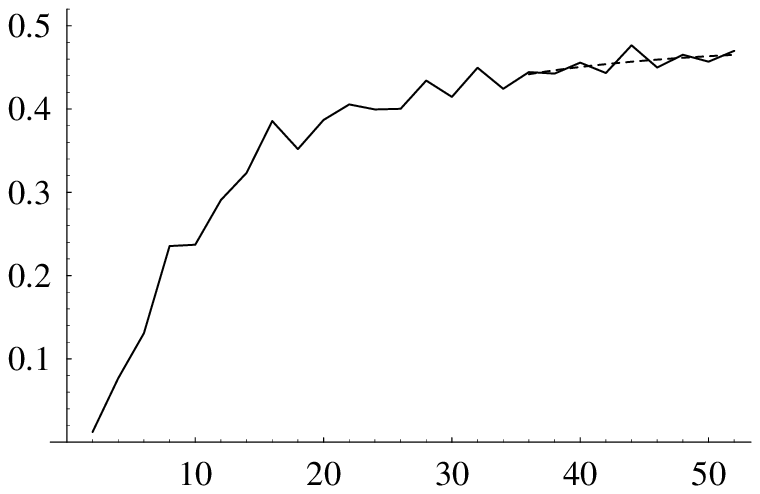}%
\par
{\footnotesize{
{\textbf{Figure 8.} $|\nu_j|^{-1/j}_3$
as a function of $j \in \{2,4,...,50,52\}$. The dashed
line is the graph of the interpolant in \rref{nu3jin}, for $j \in [36,52]$.
}}
\par}
\label{f4}
}
\end{figure}
Now, recalling that $\theta_3$ is a
lower bound on the time of existence $T$ of the solution $u$ (see \rref{claim3}), we are led to the final estimate
\beq 0.47 < T \leqs +\infty~. \feq
In particular, as anticipated, we
have indications that $u$ does not blow up near the
$H^3$-convergence radius $\tau_3$.
\vfill \eject \noindent
\section{Possible blow-up at larger times
for the Behr- Ne$\check{\mbox{c}}$as-Wu datum, via \Pade approximants}
\label{secpade}
\textbf{A few words on \Pade approximants.}
Let us be given an analytic function $f : I \vain \complessi, t \mapsto f(t)$,
with $I$ a neighborhood of zero in $\reali$ or $\complessi$.
Let $p, q \in \{0,1,2,...\}$; we recall that the \Pade approximant of order
$(p,q)$ of $f$, if it exists, is the unique complex function
$[p/q]_f \equiv [p/q]$ of the form
\beq [p/q](t) = {a_0 + a_1 t + ... + a_p t^p \over 1 + b_1 t + ... b_q t^q}~, \feq
such that
\beq f(t) = [p/q](t) + O(t^{p+q+1})~~\mbox{for $t \vain 0$}~; \feq
the above condition determines the $p+q+1$ unknown coefficients
$a_0, ..., b_q$ as functions of the derivatives $f^{(j)}(0)$, $j=0,...,p+q$;
the domain of $[p/q]$ is the largest subset of $\complessi$ where the above ratio
is defined. The family of all approximants $[p/q]$ ($p,q = 0,1,2,...$) forms
the so-called \Pade table of $f$; the approximants with $p=q$ are called diagonal. \par
There are several results and conjectures about the convergence to $f$
of the \Pade approximants $[p/q]$ with $p$ or $p,q$ large. In particular,
the so-called ``\Pade conjecture'' (or ``Baker-Gammel-Wills conjecture'')
states that, for a meromorphic function $f$ on
a disk of $\complessi$, there is a subsequence $[p_\ell/p_\ell]$ ($\ell=1,2,3,...$) of
diagonal \Pade approximants that, for $\ell \vain + \infty$, converges to $f$
uniformly on each compact subset of the disk minus the poles of $f$. This conjecture
has been proved for special classes of meromorphic functions (see
\cite{Bak} \cite{Pad1} \cite{Pad2} and references therein). \par
It is found experimentally that the \Pade approximants
of large order (and, in particular, the diagonal
approximants $[p/p]$) work as well for many
non meromorphic functions, describing accurately
their behavior even close to non polar
singularities.
\salto
\textbf{\Pade approximants for $\boma{\| u(t) \|^2_3}$, and
possible evidence for a blow-up.}
The previous considerations can be applied
(for suitable $n$) to the function $f(t) := \| u(t) \|^2_n$,
where $u$ is the solution of the Euler equation
with a given datum $u_0$. \par
One can ascribe to a number of works
the idea of using the \Pade approximants for
such a function; as in the Introduction, we mention
\cite{Bra} \cite{Fri} \cite{Fri0} \cite{Pel} (and some references therein).
As already remarked, these
papers have considered initial data $u_0$ different from
the one of Behr-Ne$\check{\mbox{c}}$as-Wu (e.g., the Taylor-Green vortex); furthermore,
they have generally considered the Sobolev norm of order $n=1$.
\par
Here we are focusing on
the (maximal $\Aa$-) solution $u$ for the Behr-Ne$\check{\mbox{c}}$as-Wu datum;
from now on, $[p/q]$ stands for the \Pade approximants of
the analytic function
\beq t \mapsto f(t) := \| u(t) \|^2_3~. \label{fu3} \feq
We conjecture that, for certain large $p$, $[p/p]$ approximates well the function
$t \mapsto \| u(t) \|^2_3$ (and even its analytic continuation
to the complex plane).
From the previous paragraphs, we have
the derivatives $f^{(j)}(0) = j!\,  \nu_{3 j}$ for $j=0,...,52$;
this information suffices
to determine all the \Pade approximants $[p/q]$ for
$p+q \leqs 52$ and, in particular, all the
diagonal approximants $[p/p]$ for $p=0,1,...,26$. \par
It turns out that the diagonal approximants $[p/p]$
exist in the cases of even order
$p=0,2,...,26$, while they do not exist
in the odd cases $p=1,3,...,25$ (the reason being,
essentially, that
the power series for $f(t)$ about zero contains only even
powers of $t$).
Let us consider, for example, the approximant $[12/12]$. Its numerator
and denominator are polynomials with rational coefficients, too
large to be written explicitly; however, we can use the $16$-digits
approximation for the coefficients and write
\parn
\vbox{
\beq [12/12](t) = {N_{12}(t) \over D_{12}(t)}~, \feq
$$ N_{12}(t) := 96 + 6.680481407149543 \times 10^3\,t^2
+ 3.08095009988031 \times 10^4 \,t^4 $$
$$ + 2.3462351635051233 \times 10^5 \,t^6
+ 2.407391215430808 \times 10^5 \,t^8 $$
$$ + 2.5575522886490226 \times 10^6\,t^{10} + 3.094974424148063 \times 10^6\,t^{12}~, $$
$$ D_{12}(t) := 1 + 0.255014657807743\,t^2 + 4.288322833232482\,t^4 - 5.985294148961588\,t^6 $$
$$ +  8.973150435320479\,t^8 + 66.29326162173366\,t^{10} - 612.1107629833056\,t^{12}~. $$
}
The poles of $[12/12]$, which are the zeros of $D_{12}$, are simple and occur at the points
\parn
\vbox{
\beq t = \pm 0.294020 \pm 0.464361 \, i \quad (|t| = 0.549617)~; \feq
$$ t = \pm 0.511609 \pm 0.301416 \,i \quad (|t| = 0.593797)~; $$
$$ t = \pm 0.606004 \, i~, \qquad t = \pm 0.626199 $$
}
(here and in the sequel, $\pm$ means that we can choose independently
the signs for the real and imaginary part, e.g.,
$+$ for the real and $-$ for the imaginary part).
So, the singularities of minimum modulus of the approximant $[12/12]$
are at anyone of the points $\to = \pm 0.294020 \pm 0.464361 \, i$,
such that $|\to| = 0.549617$; furthermore, the real singularities
closest to the origin are at anyone of the points $\ts = \pm 0.626199$. \par
We have performed a similar analysis for all the approximants
$[p/p]$, with $p=14,16,...,26$; the results are summarized in
Table 2.
\salto
\vbox{
\noindent
\textbf{Table 2.} Poles of the \Pade approximants $[p/p](t)$ to $\| u(t) \|^2_3$. \parn
$\to :=$ pole closest to the origin (with modulus $|\to|$); \parn
$\ts := $ real (or almost real) pole closest to the origin.
$$
\begin{tabular}{|l|l|l|l|}
\hline
$[p/p]$ & $\to$ & $|\to|$ & $\ts$ \\
\hline
$[12/12]$ &$\pm 0.294020 \pm 0.464361 \, i$  & $0.549617$ & $\pm 0.626199$ \\
$[14/14]$ &  $\pm 0.281333 \pm 0.445002 \, i$ & $0.526474 $ & $\pm 0.656185 $ \\
$[16/16]$ &  $ \pm 0.283300 \pm 0.446498 \, i$ & $0.528790$ & $\pm 0.661087 $ \\
$[18/18]$ &  $\pm 0.283081 \pm 0.445859 \, i$, & $0.528134$ & $\pm 0.660118 $ \\
$[20/20]$ &  $\pm 0.345307 \pm 0.348713 \, i $ & $0.490752$ & $\pm 0.621387 \pm 0.047708 \, i $ \\
$[22/22]$ &  $\pm 0.350239 \pm 0.350695  \, i $ & $0.495635$ & $\pm 0.541967 $ \\
$[24/24]$ &  $\pm 0.349063 \pm  0.350777 \, i $ & $0.494863 $ & $\pm 0.609804 \pm 0.0383530 \, i  $ \\
$[26/26]$ &  $\pm 0.0714399 \pm  0.508700 \, i $ & $0.513692 $ & $\pm  0.816133  $ \\
\hline
\end{tabular}
$$
}
\salto
Let us point out some features of the \Pade approximants Table 2,
with their possible implications: \parn
(i) For all the approximants $[p/p]$ in the table,
the poles of minimum modulus occur at
points $\to$ with $|\to| \simeq 0.5$. There is not
a clear trend of $|\to|$ as a function of $p$,
so we limit ourself to consider
the mean of $|\to|$ for $p=12,...,26$ which is
$\la |\to| \ra = 0.515995$, with a mean
quadratic error $\Delta \to < 0.02$.
On the other hand, for a holomorphic function, the convergence radius of the
power series centered at zero is the modulus of
the singularity closest to the origin. So, assuming that
the above $[p/p]$ describe approximately the singularities
of $f(t) = \| u(t) \|^2_3$, we can derive from
these approximants an estimate of a convergence
radius $\theta_3$ for the power series
of $f(t)$. More precisely, assuming
$|\to | - \Delta \to < \theta_3 < |\to |
+ \Delta \to$ and rounding up to two digits, we
obtain from the above \Pade approximants an
estimate
\beq 0.49 < \theta_3 < 0.54~; \feq
this is compatible with the estimate on $\theta_3$
obtained in Section \ref{ourap} by other means (see Eq. \rref{teta3}
and the discussion before it). \parn
(ii) The $[p/p]$ approximants of Table 2 have real poles (symmetric
with respect to the origin),
with the exceptions of $[20/20]$ and $[24,24]$ which,
however, possess ``almost real'' poles, close to the real axis
({\footnote{The occurring of almost real singularities
has also been pointed out in \cite{Fri} \cite{Pel} while analyzing
the \Pade approximants for $\| u(t) \|^2_1$, with initial conditions $u_0$ different
from the Behr-Ne$\check{\mbox{c}}$as-Wu datum.}}).
In the table, we have denoted with $\ts$ the real (or almost real) singularities
closest to the origin. The mean of $|\ts|$ for
for $p=12,...,26$ is $\la |\ts| \ra = 0.649489$, with
a mean quadratic error $\Delta \ts < 0.08$ (however,
there are large deviations from the mean in the
special cases $p=22$ and $p=26$). \parn
The above results on the singularities $\ts$ somehow suggest
that $f(t) = \| u(t) \|^2_3$ could diverge for $t \vain T^{-}$ (and $t
\vain (- T)^{+}$), for a suitable $T$; if we assume for $T$
the upper and lower bounds $|\ts| \pm \Delta \ts$, rounding
up to two digits we get
\beq 0.56 < T < 0.73~. \label{boundt} \feq
If such a conjectured divergence of $f(t)$ actually occurred, the solution $u$ of the Euler equation with the
Behr-Ne$\check{\mbox{c}}$as-Wu datum would blow up at $T$ (and $-T$);
admittedly, the indications for such a blow up are very weak.
\salto
\textbf{D-log Pad\`e approximants for $\boma{\| u(t) \|^2_3}$}. As well known,
the D-log \Pade approximants of a function $t \mapsto f(t)$ are
the \Pade approximants for the logarithmic derivative $\dot{f}/f$
($\dot{~} := d/dt$).
These approximants are generally regarded as more
suitable for describing the behavior of $f$ close
to singularities, even of non polar type. In particular, the
presence of a singularity at a point $T_{*}$, say
real, and a behavior of the type
$[p/p]_{\dot f/f} \sim \lambda_{*}/(T_{*} - t)$
for $t \vain T^{-}_{*}$
is regarded as an indication that $f(t) \sim \mbox{const}/(T - t)^{\lambda}$
for real $t \vain T^{-}$, where $T \simeq T_{*}$ and $\lambda \simeq \lambda_{*}$ \cite{Bak}. \par
We have attempted an analysis of the function $f(t) := \| u(t) \|^2_3$
via the approximants $[p/p]_{\dot f/f}$, with odd $p \leqs 25$
({\footnote{Our function has the form $f(t) = F(t^2)$;
in such a case, for odd $p$, the D-log approximant of $f$ of
order $(p,p)$ is (up to a factor $2 t$)
the D-log approximant of the function $s \mapsto F(s)$ of order
$(p/2-1/2, p/2-1/2)$. On the contrary, for even $p$,
the $(p,p)$ D-log approximant of $f$ cannot be interpreted in terms
of $F$.}}); the results are very unstable with respect to the order,
and ultimately not sufficient to get any indication of blow-up.
(\footnote{Here is a more precise description
of the computational outcomes.
The D-log approximants of order $(p,p)$ for $p=17,19,21$ have
real singularities at points $T_{*} \simeq 0.72$
and are such that $[p/p]_{\dot f/f} \sim \lambda_{*}/(T_{*} - t)$
for $t \mapsto T^{-}_{*}$, with $\lambda_{*} \simeq 2.6$;
so, for $\| u(t) \|_3 = \sqrt{f(t)}$ we have a conjecture
$\| u(t) \|_3  \sim \mbox{const.}/(T - t)^{\alpha}$
with $T \simeq 0.72$ and $\alpha \simeq \lambda_{*}/2 \simeq 1.3$.
This value of $\alpha$ agrees with the Beale-Kato-Majda bound $\alpha \geqs 1$ in
the event of
blow-up (see Eq. \rref{teorkato}); it agrees as well with the (conjectural)
bound $\alpha \geqs 6/5$, obtained extrapolating from $\reali^3$ to $\Tt$
the estimate \rref{teorest}.\par
On the contrary, the D-log approximant of order $(23,23)$ for $f$ has
no real (nor almost real) singularity. Finally, at the order $(25,25)$ there is a real
singularity for $T_{*} \simeq 0.52$, and
$[25/25]_{\dot f/f} \sim \lambda_{*}/(T_{*} - t)$
for $t \mapsto T_{*}^{-}$, with $- 0.002 < \lambda_{*} < 0.002$
(there are numerical difficulties in a more precise determination
of $\lambda_{*}$). Returning to $\| u(t) \|_3 = \sqrt{f(t)}$,
the $[25, 25]$ \Pade would suggest
$\| u(t) \|_3  \sim \mbox{const.}/(T - t)^{\alpha}$
with $T \simeq 0.52$ and $-0.001 \lesssim \alpha \lesssim 0.001$.
This statement is an absurdity even in the case
$0 < \alpha \lesssim 0.001$, since it contradicts
the Beale-Kato-Majda bound \rref{teorkato} $\alpha \geqs 1$.})
\vfill \eject \noindent
\section{Conclusions}
The previous results about
the Behr-Ne$\check{\mbox{c}}$as-Wu datum $u_0$ support
our statements in the Introduction, i.e.: \parn
\vbox{
\noindent
(a) The power series for $u_0$ has an $H^3$
convergence radius $\tau_3$ such that $0.32 < \tau_3  < 0.33$
(see Eq. \rref{estau3}). \parn
(b) There is no blow-up
at time $\tau_3$ and the (maximal $\Aa$-) solution $u$ of the Euler
Cauchy problem exists, at least, up
to a time $\theta_3$ (the convergence
radius for the series expansion of $\| u(t) \|^2_3$),
for which we have from \rref{teta3} the estimate $\theta_3 > 0.47$. \parn
(c) The \Pade approximants for $\| u(t) \|^2_3$ in Table 2 give weak indications
that $u$ might blow up at a time $T$, with $0.56 < T < 0.73$ (see Eq.
\rref{boundt}).}
\parn
We think that the evidence given in this paper is rather strong for (a)(b).
As for (c), doubts on the blow-up conjecture arise not only from
the rather erratic behavior of the real singularities in the
computed \Pade approximants; in fact there are more general reasons,
recalled at the end of the Introduction, suggesting caution in
deriving blow-up results from the \Pade approximants.
\vskip 0.7cm \noindent
\textbf{Acknowledgements.}
We are grateful to Paolo Butera for useful bibliographical indications
and appreciated expertise on \Pade approximants.
\parn
This work was supported by INdAM, INFN and by MIUR, PRIN 2008
Research Project "Geometrical methods in the theory of nonlinear waves and applications".

\end{document}